\documentclass[11pt]{article}

\usepackage{natbib}
\usepackage{latexsym}
\usepackage{amssymb}
\usepackage{amsmath}
\usepackage{graphicx}
\usepackage{comment}
\usepackage{psfrag}
\usepackage{float}
\usepackage{color}
\usepackage{nicefrac}
\usepackage{ntheorem}
\usepackage{hyperref}
\usepackage{tikz}
\usetikzlibrary{backgrounds}

\def\nr{\par }

\def\beq{\begin{equation}}
\def\eeq{\end{equation}}

\newtheorem{theorem}{Theorem}
\newtheorem{lemma}{Lemma}
\newtheorem{corollary}{Corollary}

\newtheorem{proposition}{Proposition}
\newtheorem{assumption}{Assumption}
\newtheorem{definition}{Definition}

\newtheorem{example}{Example}
\newtheorem{remark}{Remark}
\newcommand{\proof}{\bf Proof: \rm \nr}

\newcommand{\qed}{\hfill $\Box$ \nr \medskip}

\def\ba{\begin{array}}
\def\ea{\end{array}}
\def\beann{\begin{eqnarray*}}
\def\eeann{\end{eqnarray*}}
\def\bea{\begin{eqnarray}}
\def\eea{\end{eqnarray}}

\def\BT{\begin{theorem}}
\def\ET{\end{theorem}}
\def\BL{\begin{lemma}}
\def\EL{\end{lemma}}
\def\BC{\begin{corollary}}
\def\EC{\end{corollary}}
\def\BE{\begin{example}}
\def\EE{\end{example}}
\def\BD{\begin{definition}}
\def\ED{\end{definition}}
\def\BR{\begin{remark}}
\def\ER{\end{remark}}
\def\BAS{\begin{assumption}}
\def\EAS{\end{assumption}}
\def\BI{\begin{itemize}}
\def\EI{\end{itemize}}

\def\BMP{\begin{minipage}{9.5cm}}
\def\EMP{\end{minipage}}
\def\MPT{\begin{minipage}{11.5cm}}
\def\EPT{\end{minipage}}

\def\R{\mathbb{R}}

\newcommand*\samethanks[1][\value{footnote}]{\footnotemark[#1]}

\begin{document}
\title{On stability of the Scholtes regularization for\\  mathematical programs with complementarity constraints}
\author{
S. L\"ammel
\thanks{
Department of Mathematics, Chemnitz University of Technology,
Reichenhainer Str. 41, 09126
Chemnitz, Germany; e-mail: sebastian.laemmel@mathematik.tu-chemnitz.de, vladimir.shikhman@mathematik.tu-chemnitz.de.
 } \and V. Shikhman\samethanks[1]
}

\maketitle

\maketitle
\vspace{-5ex}
\abstract{
%For the class of mathematical programs with complementarity constraints (MPCC), we refine the convergence analysis of the Scholtes regularization known from the literature. Our goal is to relate nondegenerate C-stationary points of MPCC with nondegenerate Karush-Kuhn-Tucker points of the corresponding Scholtes regularization. By investigating this issue, we detect the following anomalies: (i) in an arbitrarily small neighborhood of a nondegenerate C-stationary point there could be degenerate Karush-Kuhn-Tucker points of the Scholtes regularization; (ii) even if nondegenerate, they might though be locally non-unique; (iii) if nevertheless unique, their quadratic index potentially differs from the C-index of the C-stationary point under consideration. As a striking consequence of this index shift, a change of the topological type for Karush-Kuhn-Tucker points of the Scholtes regularization is possible. In particular, a nondegenerate minimizer of MPCC might well be approximated by the saddle points of its Scholtes regularization. In order to bypass the mentioned anomalies, an additional generic condition for nondegenerate C-stationary points of MPCC is identified. By assuming the latter, we uniquely trace then nondegenerate Karush-Kuhn-Tucker points of the Scholtes regularization and successively maintain their topological type.

For mathematical programs with complementarity constraints (MPCC), we {\color{black} study the stability properties} of their Scholtes regularization. Our goal is to relate nondegenerate C-stationary points of MPCC with nondegenerate Karush-Kuhn-Tucker points of the Scholtes regularization {\color{black} up to their topological type. 
As it is standard in the framework of Morse theory, the topological types are captured by the C-index and the quadratic index, respectively. } %Index is the main invariant in optimization in order to structurally distinguish local minimizers / maximizers and saddle points from each other. Our first observation reveals the so-called index shift phenomenon. It may well happen that the C-index of a C-stationary point of MPCC differs from the quadratic index of the approximating Karush-Kuhn-Tucker points of the Scholtes regularozation. 
 {\color{black} It turns out that}
 a change of the topological type for {\color{black} the approximating }Karush-Kuhn-Tucker points of the Scholtes regularization {\color{black} and their limiting C-stationary point is possible.} In particular, a minimizer of MPCC {\color{black} with zero C-index} might be approximated by saddle points {\color{black} of the Scholtes regularization with nonzero quadratic index}. In order to bypass this {\color{black} index shift phenomenon}, an additional generic condition for nondegenerate C-stationary points of MPCC is identified. {\color{black} It says that non-biactive multipliers under consideration should not vanish.} Then, we uniquely trace nondegenerate Karush-Kuhn-Tucker points of the Scholtes regularization and successively maintain the topological type of their limiting C-stationary point. {\color{black} The main technical issue here is to relate the first-order information of the defining functions, which enters the biactive part of the C-index, with the second-order information, which enters the quadratic index of the Karush-Kuhn-Tucker points.}
%As a byproduct, we refute a result on the well-posedness of Scholtes regularization wrongly proven in \citet{still:2007}.
}%

\vspace{2ex}

\vspace{2ex}
{\bf MSC-classification: 90C26, 90C33
}

% Text of your paper here
\section{Introduction}\label{sec:intro}
%The seminal paper \cite{scheel:2000} started a line of research. In several papers e.g~ \cite{still:2007}, \cite{ralph:2004}

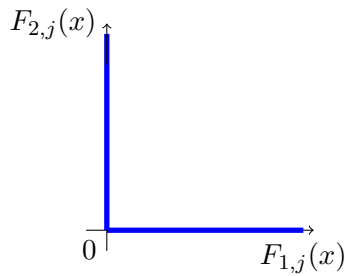
\begin{figure}[h]
	\centering
	\begin{tikzpicture}[xscale=0.55,yscale=0.55,domain=-0.5:5,samples=100]
	\draw node[below left] {0};
	\draw[->] (-0.5,0) -- (5,0);
	\draw[blue,line width=2] (0,0) -- (4.75,0) node[below,black] {$F_{1,j}(x)$};
	\draw[blue, line width=2] (0,0) -- (0,4.75) ;
	\draw[] (0,-0.5) -- (0,0);
	\draw[->] (0,4) -- (0,5) node[left] {$F_{2,j}(x)$};
	\end{tikzpicture}
	\caption{MPCC feasible set $M$}  
	\label{fig:feasiblesetMPCC}
\end{figure}

Mathematical programs with complementarity constraints (MPCC) have been widely studied in the literature for several decades, see e.g.~\cite{pang:2003}:
{\color{black}
%We consider mathematical programs with complementarity constraints:
\[
\mbox{MPCC}: \quad
\min_{x} \,\, f(x)\quad \mbox{s.\,t.} \quad x \in M
\]
with
\[
M=\left\{x \in\R^n\, \left\vert\,\
F_{1,j}(x) \cdot F_{2,j}(x)=0, F_{1,j}(x) \ge 0, F_{2,j}(x)\ge 0, j=1,\ldots,\kappa \right.\right\},
\]
where the defining functions $f \in C^2(\R^n,\R)$,
$F_1, F_2 \in C^2(\R^n,\R^\kappa)$ are twice continuously differentiable. For simplicity, we omit additional smooth equality and inequality constraints here, since a generalization of our results to this case is straightforward.
%We refer to $f, F_1$, and $F_2$ as the defining functions. If needed, we highlight the dependence on them by $\mbox{MPCC}(f,F_1,F_2)$.
}
MPCCs naturally appear if dealing with 
generalized Nash equilibrium problems, bilevel optimization, and semi-infinite programming to name a few applications. The main difficulty of dealing with MPCCs is that their feasible sets {\color{black} exhibit nonsmoothness}. {\color{black} A typical MPCC feasible set $M$ is depicted in Figure \ref{fig:feasiblesetMPCC}, where $n=2$, $\kappa=1$, $F_{1,1}(x)=x_1$ and $F_{2,1}(x)=x_2$.} 

%It is well-known that the MPCC feasible set can be equivalently described by min-type equality constraints. 
%However, in general most of the standard constraints qualification fail to hold at any feasible point.

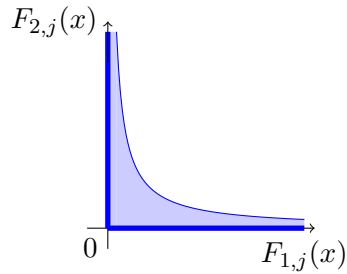
\begin{figure}[h]
    \centering
\begin{tikzpicture}[xscale=0.55,yscale=0.55,domain=-0.5:5,samples=100]
    \draw node[below left] {0};
    \fill [blue!20]
        (0, 0)
      -- plot[domain=0:1/4.75] (\x, {4.75})
      -- (1/4.75, 0);
    \fill [blue!20]
        (1/4.75, 0)
      -- plot[domain=1/4.75:4.75] (\x,{1/ \x})
      -- (4.75, 0);
    \draw[blue] plot[domain=1/4.75:4.75] (\x,{1/ \x});
    \draw[->] (-0.5,0) -- (5,0);
        \draw[->] (0,4) -- (0,5) node[left] {$F_{2,j}(x)$};
     \draw[blue,line width=2] (0,0) -- (4.75,0) node[below,black] {$F_{1,j}(x)$};
    \draw[blue, line width=2] (0,0) -- (0,4.75) ;
    \draw[] (0,-0.5) -- (0,0);
\end{tikzpicture}
\caption{Feasible set $M^\mathcal{S}$ of the Scholtes regularization}    \label{fig:feasiblesetS}
\end{figure}

The nonsmoothness of the MPCC feasible set motivates the introduction of regularization methods, see e.g.~\cite{hoheisel:2013}.
In the seminal paper \cite{scholtes:2001}, such a regularization due to Scholtes was suggested first:
{\color{black}
\[
\color{black}
\mathcal{S}:
\color{black}\quad
\min_{x} \,\, f(x)\quad \mbox{s.\,t.} \quad x \in M^\mathcal{S}
\]
with
\[
    M^\mathcal{S}=\left\{x \in\R^n\, \left\vert\,\begin{array}{l}
    F_{1,j}(x) \cdot F_{2,j}(x)\le t, F_{1,j}(x) \ge 0,F_{2,j}(x)\ge 0, j=1,\ldots,\kappa
    \end{array} \right.\right\},
\]
where the parameter $t>0$ is positive and is suppressed in the notation for the readers' convenience. In Figure \ref{fig:feasiblesetS}, the Scholtes regularization of the typical MPCC feasible set is depicted.
}
Under suitable assumptions, two types of results -- on its convergence and well-posedness -- have been proven there and further in \cite{ralph:2004}: (1) the Karush-Kuhn-Tucker points of the Scholtes regularization converge to a C-stationary point of MPCC; (2) such a sequence of Karush-Kuhn-Tucker points exists in a neighborhood of a C-stationary point.
Up to this date numerical analyses have shown the superiority of the Scholtes regularization
even over more modern and theoretically better regularizations \cite{hoheisel:2013, kanzow:Epsilon, nurkanovic:2024}.

In this paper, we refine the convergence analysis of the Scholtes regularization performed in the literature so far. 
Our goal is to relate nondegenerate C-stationary points of MPCC with nondegenerate Karush-Kuhn-Tucker points of the corresponding Scholtes regularization up to their topological type. To put it simply,  nondegeneracy refers to the validity of tailored constraint qualifications, non-vanishing of certain multipliers, and second-order regularity. In the framework of nonlinear programming (NLP), the topological type of a nondegenerate Karush-Kuhn-Tucker point is characterized by its quadratic index, i.e.~the number of negative eigenvalues of the restricted Hessian matrix of the Lagrange function, see e.g.~\cite{jongen:2000}. 
For a nondegenerate C-stationary point of MPCC, the number of negative pairs of biactive multipliers is additionally involved. This biactive part, together with the quadratic part as before, yields the C-index. Again, it can be unambiguously deduced from the C-index if the corresponding C-stationary point  is a minimizer, maximizer, or a saddle of a certain type. As shown e.g.~in \cite{shikhman:2012}, the C-index for MPCC is the right generalization of the quadratic index for NLP from the topological point of view. Namely, it adequately describes the change of the MPCC lower level sets up to homotopy equivalence in the sense of Morse theory. 

 Let us briefly motivate our studies. While considering the Scholtes regularization of an MPCC instance, it is not only important to relate minimizers of both. Also invariants of the corresponding stationary points, such as their nondegeneracy, uniqueness, and index, become crucial. This is in order to address the global structure of the Scholtes regularization for MPCC and, finally, to identify conditions under which MPCCs can be regularized in a stable manner.
{\color{black} The direct attempt to such a study reveals certain instabilities:} 
%Our topological approach to the study of convergence and well-posedness properties of the Scholtes regularization reveals several anomalies: 
\begin{itemize}
%    \item [(i)] Nondegeneracy of a sequence of Karush-Kuhn-Tucker points is not necessarily preserved in the limiting point, i.e.~the obtained C-stationary point of MPCC might be degenerate, see Example \ref{ex:ndc2fail}. Vice versa, it can not be assured to find nondegenerate Karush-Kuhn-Tucker points of the Scholtes regularization in any arbitrarily small neighborhood of a nondegenerate C-stationary point of MPCC, see Example \ref{ex:ssosc}.
    \item [(i)] The Scholtes regularization might cause bifurcation, i.e.~in a neighborhood of a {\color{black}nondegenerate} C-stationary point of MPCC
    there are potentially multiple Karush-Kuhn-Tucker points of the former{\color{black}, see Example \ref{ex:2min}}. This includes cases where they form a continuum, see Example \ref{ex:continuum}, or the C-stationary point under consideration is a minimizer of MPCC, see Examples \ref{ex:continuum2}. %Thus the structure of the regularization might become more complicated than that of the underlying MPCC.
    \item [(ii)] The topological type might change while regularizing. {\color{black} This is due to the so-called index shift phenomenon.} The latter means that the C-index of the limiting C-stationary point might decrease if compared to the quadratic index of the approximating Karush-Kuhn-Tucker points. We refer to Example \ref{ex:NDC4}, where saddle points of the Scholtes regularization converge to a nondegenerate minimizer of MPCC. 
\end{itemize}

In order to avoid {\color{black} both, bifurcation and index shift}, an additional generic condition, called by us NDC4, for nondegenerate C-stationary points of MPCC is identified. It requires that the multipliers associated with non-biactive constraints do not vanish either.
By additionally assuming this, we uniquely trace then nondegenerate Karush-Kuhn-Tucker points of the Scholtes regularization and successively maintain their topological type in comparison with the limiting C-stationary point of MPCC. For the corresponding convergence result see Theorem \ref{thm:kktsequence} and for the corresponding well-posedness result see Theorem \ref{thm:wellposedness}. {\color{black} The challenge behind the proofs of these results is to relate the quadratic index of Karush-Kuhn-Tucker points of the Scholtes regularization with the C-index of C-stationary points of MPCC, and vice versa. The main difficulty here is that the biactive part of the latter is given by the first-order, whereas the quadratic index of the former by the second-order information of the defining functions.}

{\color{black}
In what follows, we provide intuition about NDC4, which by definition prevents non-biactive multipliers from vanishing. We  emphasize that this additional condition quite naturally arises in the context of Scholtes regularization. For seeing this, it is enough to set $\kappa=1$ and to consider just one complementarity constraint
\begin{equation}
    \label{eq:intui}
  F_{1,j}(x) \cdot F_{2,j}(x) =0, \quad F_{1,j}(x) \geq 0, \quad F_{2,j}(x) \geq 0,
\end{equation}
along with its Scholtes regularization
\begin{equation}
    \label{eq:intui-t}
  F_{1,j}(x) \cdot F_{2,j}(x) \leq t, \quad F_{1,j}(x) \geq 0, \quad F_{2,j}(x) \geq 0.
\end{equation}
Let us  assume that at a feasible point $\bar x \in M$ satisfying (\ref{eq:intui}) we have $F_{1,j}(\bar x)=0$ and $F_{2,j}(\bar x) >0$, without loss of generality. Then, locally around $\bar x$, $F_{1,j}(x) \cdot F_{2,j}(x) =0$ in (\ref{eq:intui}) reduces to the equality constraint $F_{1,j}(x)=0$, and the constraint $F_{1,j}(x)\geq 0$ becomes redundant. The situation differs if we let $t \rightarrow 0$ in (2):
\begin{equation}
    \label{eq:intui-0}
  F_{1,j}(x) \cdot F_{2,j}(x) \leq 0, \quad F_{1,j}(x) \geq 0, \quad F_{2,j}(x) \geq 0.
\end{equation}
Clearly, both descriptions (\ref{eq:intui}) and (\ref{eq:intui-0}) yield the same feasible set. However, $F_{1,j}(x) \cdot F_{2,j}(x) \leq 0$ in (\ref{eq:intui-0}) now locally reduces to the inequality constraint $F_{1,j}(x)\leq 0$, and the absence of the constraint $F_{1,j}(x)\geq 0$ would lead to a different feasible set. Therefore, it makes sense that also the multiplier information of the non-biactive constraint becomes important for the limiting description (\ref{eq:intui-0}). %Note that it is however superfluous for the original description (\ref{eq:intui}).
}
{\color{black} To elaborate more on this, let us additionally assume that $\bar x \in M$ is a nondegenerate C-stationary point of MPCC with the non-biactive multiplier $\bar \sigma_{1,j}$ corresponding to the locally active constraint $F_{1,j}(x)=0$. First, we focus on the case $\bar \sigma_{1,j} >0$. Then, $\bar x$ is a nondegenerate Karush-Kuhn-Tucker point with respect to the feasible set given just by the inequality constraint $F_{1,j}(x)\geq 0$. Hence, $\bar x$ is also a nondegenerate Karush-Kuhn-Tucker point of the Scholtes regularization for any $t > 0$, cf. Figure \ref{fig:feasiblesetS}. We point out that NDC4 here, i.e. $\bar \sigma_{1,j} \not = 0$, corresponds to the classical notion of nondegeneracy for NLP, i.e. multipliers for the active inequality constraints should not vanish, see e.g. \cite{jongen:2000}. Second, we consider the case $\bar \sigma_{1,j} <0$. If at all, $\bar x$ must then be approximated by those Karush-Kuhn-Tucker points $x^t \in M^\mathcal{S}$ of the Scholtes regularization for $t \rightarrow 0$, which fulfill $F_{1,j}(x^t) \cdot F_{2,j}(x^t)=t$. This is since NDC4 at $\bar x$ rules out the approximating Karush-Kuhn-Tucker points of the Scholtes regularization lying in the interior of $M^\mathcal{S}$, cf. Example \ref{ex:NDC4}. 
}
%%%%%%%%%%%%%%%%%%%%%%%%%% TO DO

The obtained results on the Scholtes regularization of MPCC perfectly fit into the existing literature and considerably improve on the state of art. For the well-posedness part, our Theorem \ref{thm:wellposedness} generalizes the results in \cite{scholtes:2001} from the case of minimizers to all types of nondegenerate C-stationary points of MPCC. 
%Moreover, the bifurcation phenomenon  in the absence of NDC4 enables us to refute the incorrect uniqueness part in \cite[Theorem 5.1]{still:2007}, see Example \ref{ex:2min}.  { More than that, we emphasize that the proof  of existence there also fails. This has a negative influence on more recent research on related problems, see \cite{bouza:2013} for bilevel programming and \cite{bouza:2023} for one-parametric MPCCs, where \cite[Theorem 5.1]{still:2007} has been used.} 
Moreover, we make a comparison with the well-posedness result of \cite{ralph:2004}, which also generalizes that of \cite{scholtes:2001}. It turns out that the weaker assumptions from \cite{ralph:2004} do not prevent from degeneracies and, hence, are not suitable for the aims of our studies, see Example \ref{ex:ssosc}.
The 
%of the "somewhat artificial" (as stated in \cite{ralph:2004} 
conditions of nonvanishing multipliers, including NDC4, are therefore not "somewhat artificial", as claimed in \cite{ralph:2004}, but rather essential for stability of the Scholtes regularization.
In addition, we emphasize that our convergence results in Theorem \ref{thm:kktsequence} are novel due to the discovery of the index shift phenomenon.
 {
They suggest that it is not sufficient to focus on biactive multipliers as done e.g. in \cite{wang:2023} with a numerical study of MPCC, but on non-biactive multipliers as well.}
{
\color{black}
Finally, we apply our analysis to the smoothing approach from \cite{still:2007}, where inequality constraints in (\ref{eq:intui-t})
are substituted by equality constraints:
\begin{equation}
    \label{eq:intui-eq}
  F_{1,j}(x) \cdot F_{2,j}(x) = t, \quad F_{1,j}(x) \geq 0, \quad F_{2,j}(x) \geq 0.
\end{equation}
In this setting, NDC4 is not required anymore for uniquely tracing the topological type of the corresponding C-stationary and Karush-Kuhn-Tucker points, see Theorems \ref{thm:kktstill} and \ref{thm:wpstill}. As above, this can be explained by taking $t \rightarrow 0$ in (\ref{eq:intui-eq}). Since  obtaining the original description (1) in the limit, the usual nondegeneracy conditions for MPCC suffice for guaranteeing stability of the regularization (\ref{eq:intui-eq}). We conclude that NDC4 is rather tailored for the description of the MPCC feasible set as given in (\ref{eq:intui-0}) and is therefore crucial for dealing with stability aspects of the Scholtes regularization.
}

%%%%%%%%%%%%

The article consists of four sections. We discuss in Section \ref{sec:preliminaries} preliminaries for MPCC and its Scholtes regularization. In Section \ref{sec:main}, we present our findings on the convergence behavior of the Scholtes regularization and its well-posedness properties. Subsequently, Section \ref{sec:compare} is devoted to a comparison with the existing literature.
Our notation is standard. We denote the $n$-dimensional Euclidean space by $\R^n$ with the $j$-th coordinate vector denoted by $e_j$, $j=1,\ldots,n$. Given a twice continuously differentiable function $f:\R^n\to \R$, \textcolor{black}{we denote its differential by $D f$} and its Hessian matrix by $D^2f$. The sign of $x \in \R$ is denoted by $\mbox{sgn}(x)$.

%%%%%%%%%%%%%%%%%%%%%%%%%%%%%%%%%%%%%%
%Preliminaries
%%%%%%%%%%%%%%%%%%%%%%%%%%%%%%%%%%%%%%

\section{Preliminaries}
\label{sec:preliminaries}

\subsection{Preliminaries on MPCC}\label{sec:pre}

%%%%%%%%%%%%%%%%%%%%%%%%%%%%%%%%%%%%%%
%C-stat
%%%%%%%%%%%%%%%%%%%%%%%%%%%%%%%%%%%%%%

In what follows, we briefly recall some basic notions from MPCC theory, see e.g.~\cite{shikhman:2012}. 
For that, the index sets of active constraints associated with $\bar x \in M$ will be helpful:
\[
\begin{array}{c}
a_{01}\left(\bar x\right)=\left\{j\,\left\vert\,  F_{1,j}(\bar x)=0, F_{2,j}(\bar 
x)>0\right.\right\},\\
a_{10}\left(\bar x\right)=\left\{j\,\left\vert\,  F_{1,j}(\bar x)>0, F_{2,j}(\bar  x)=0\right.\right\}, \\
a_{00}\left(\bar x\right)=\left\{j\,\left\vert\,  F_{1,j}(\bar x)=0, F_{2,j}(\bar  x)=0\right.\right\}.  
\end{array}
\]

We start our exposition with the MPCC-tailored linear independence constraint qualification.

\begin{definition}[MPCC-LICQ]
We say that a feasible point $\bar x \in M$ satisfies the MPCC-tailored linear independence constraint qualification (MPCC-LICQ) if the following vectors are linearly independent:
\[
D F_{1,j}\left(\bar x\right), j \in a_{01}\left(\bar x\right)\cup a_{00}\left(\bar x\right),\quad
D F_{2,j}\left(\bar x\right), j \in a_{10}\left(\bar x\right)\cup a_{00}\left(\bar x\right).
\]
\end{definition}
It is well-known that MPCC-LICQ is not restrictive in the sense of genericity.  
The subset of defining functions $F_1, F_2$, for which MPCC-LICQ is fulfilled at all
feasible points $\bar x \in M$ of a corresponding MPCC,  is $C^2_s$-open and -dense with respect to the strong (or Whitney-) topology, see \cite{scholtes:2001a}.

Next, we proceed with the different stationarity notions for MPCC used in the literature.

\begin{definition}[Stationarity]
    \label{def:c-stat}
A feasible point $\bar x \in M$ is called W-stationary for MPCC if there exist multipliers
\[
\begin{array}{l}
\bar \sigma_{1,j}, j \in a_{01}\left(\bar x\right),\quad \bar \sigma_{2,j}, j \in a_{10}\left(\bar x\right),\quad \bar \varrho_{1,j},\bar \varrho_{2,j}, j \in a_{00}\left(\bar x\right),
\end{array}
\]such that 
\begin{equation}
   \label{eq:cstat-1} 
   \begin{array}{rcl}
   \displaystyle D f(\bar x)&=& \displaystyle \sum\limits_{j \in a_{01}\left(\bar x\right)} \bar \sigma_{1,j} D F_{1,j}\left(\bar x\right)
    + \sum\limits_{j \in a_{10}\left(\bar x\right)} \bar \sigma_{2,j}D  F_{2,j}\left(\bar x\right) \\ \\
    &&+ \displaystyle\sum\limits_{j \in a_{00}\left(\bar x\right)} \left(\bar \varrho_{1,j} D F_{1,j}\left(\bar x\right)+\bar \varrho_{2,j} D F_{2,j}\left(\bar x\right)\right). \end{array}
\end{equation}
Moreover, the W-stationary point $\bar x$ is called:
\begin{itemize}
    \item C-stationary if
\begin{equation}
   \label{eq:cstat-2} \bar \varrho_{1,j} \cdot \bar \varrho_{2,j}\ge 0 \mbox{ for all } j \in a_{00}\left(\bar x\right);
\end{equation}
\item M-stationary if 
\begin{equation}
       \label{eq:mstat-2} \bar \varrho_{1,j} > 0, \bar \varrho_{2,j} > 0 \mbox{ or } \bar \varrho_{1,j} \cdot \bar \varrho_{2,j} = 0\mbox{ for all } j \in a_{00}\left(\bar x\right);
\end{equation}
    \item S-stationary if 
\begin{equation}
       \label{eq:sstat-2} \bar \varrho_{1,j} \ge 0, \bar \varrho_{2,j}\ge 0 \mbox{ for all } j \in a_{00}\left(\bar x\right).
\end{equation}
\end{itemize}
\end{definition}
It is clear that these stationarity notions are related as follows:
\[
\mbox{S-stationarity} \Rightarrow \mbox{M-stationarity} \Rightarrow \mbox{C-stationarity} \Rightarrow \mbox{W-stationarity}. 
\]
Moreover, S-stationarity turns out to be necessary for optimality. If $\bar x \in M$ is a local minimizer of MPCC satisfying MPCC-LICQ, then it is S-stationary, see \cite{scheel:2000}.
For a W-stationary point $\bar x \in M$ with multipliers $(\bar \sigma, \bar \varrho)$ -- which are unique under MPCC-LICQ -- it is convenient to define the appropriate Lagrange function as
\[
 \begin{array}{rcl}
\displaystyle L(x)&=&\displaystyle f(x)-\sum\limits_{j \in a_{01}\left(\bar x\right)} \bar \sigma_{1,j} F_{1,j}\left( x\right)
    - \sum\limits_{j \in a_{10}\left(\bar x\right)} \bar \sigma_{2,j} F_{2,j}\left( x\right) \\ \\
    &&- \displaystyle\sum\limits_{j \in a_{00}\left(\bar x\right)} \left(\bar \varrho_{1,j}  F_{1,j}\left(x\right)+\bar \varrho_{2,j} F_{2,j}\left(x\right)\right).
\end{array}
\]
%
\begin{comment}
  Further, we set
\[
  M_0(\bar x) =\left\{
x \in \R^n\,\left\vert\,\begin{array}{l}
F_{1,j}(\bar x)=0, j \in a_{01}\left(\bar x\right)\cup a_{00}\left(\bar x\right),\\
F_{2,j}(\bar x)=0, j \in a_{10}\left(\bar x\right)\cup a_{00}\left(\bar x\right)
\end{array} 
\right.\right\}.
\]  
\end{comment}
%
The corresponding tangent space is given by
\[
\mathcal{T}_{\bar x}=\left\{\xi \in \R^n \,\left\vert \,
\begin{array}{l}
DF_{1,j}\left(\bar x\right)\xi=0, j \in a_{01}\left(\bar x\right)\cup a_{00}\left(\bar x\right),\\
DF_{2,j}\left(\bar x\right)\xi=0, j \in a_{10}\left(\bar x\right)\cup a_{00}\left(\bar x \right)
\end{array}
\right.\right\}.
\]

Now, we turn our attention to the nondegeneracy notion for C-stationary points.

\begin{definition}[Nondegenerate C-stationarity]
    A C-stationary point $\bar x \in M$ of MPCC with multipliers $(\bar \sigma, \bar \varrho)$ is called nondegenerate if

NDC1: MPCC-LICQ holds at $\bar x$,

NDC2: multipliers for biactive constraints do not vanish, i.e. $\bar\varrho_{1,j}\cdot \bar\varrho_{2,j}> 0$,\,$j\in a_{00}\left(\bar x\right)$,
    
NDC3: the restricted Hessian matrix $D^2 L(\bar x)\restriction_{\mathcal{T}_{\bar x}}$ is nonsingular.

\noindent
For a nondegenerate C-stationary point we use the following additional condition:

NDC4: $\bar \sigma_{1,j}\ne 0$ for all $j \in a_{01}\left(\bar x\right)$ and $\bar \sigma_{2,j}\ne 0$ for all $j\in a_{10}\left(\bar x\right)$.
\end{definition}
%
% ZUM BEGRIFF ND zitieren
%We note that for a C-stationary point $\bar x$ of MPCC, NDC2 and ULSC are equivalent.
%
We note that nondegeneracy of C-stationary points is also generic.  
The subset of MPCC defining functions $f, F_1, F_2$, for which each
C-stationary point of a corresponding MPCC is nondegenerate, i.e.~fulfills NDC1-NDC3, is $C^2_s$-open and -dense, see \cite{jongen:2009}. 
For a nondegenerate C-stationary point its C-index becomes a crucial invariant, cf. \cite{ralph:2011}.

\begin{definition}[C-index]
\label{def:c-index}
    Let $\bar x \in M$ be a nondegenerate C-stationary point of MPCC with unique multipliers $\left(\bar \sigma,\bar \varrho\right)$. The number of negative eigenvalues of the matrix $D^2 L(\bar x)\restriction_{\mathcal{T}_{\bar x}}$ is called its quadratic index ($QI$). The number of negative pairs $\bar \varrho_{1,j}, \bar \varrho_{2,j}, j\in a_{00}\left(\bar x\right)$ with $\bar \varrho_{1,j}, \bar \varrho_{2,j} <0 $ is called the biactive index ($BI$) of $\bar x$. We define the C-index ($CI$) as the sum of both, i.\,e. $CI=QI+BI$.
\end{definition}
We emphasize that C-stationary points are topologically relevant in the sense of Morse theory, see \cite{jongen:2009}. It is to say that they adequately describe the topological changes of lower level sets of MPCC. First, outside the set of C-stationary points, the topology of the MPCC lower level sets remains unchanged, i.e.~the larger level set is homeomorphic to the smaller one. \textcolor{black}{Second, if passing a nondegenerate C-stationary point, a cell is attached along the boundary of the smaller lower level set to get the larger one up to the homotopy equivalence.} The dimension of the cell to be attached corresponds to the C-index of the nondegenerate C-stationary point under consideration. \textcolor{black}{Additionally, the C-index encodes also the local structure of MPCC in the vicinity of nondegenerate C-stationary points.} Nondegenerate C-stationary points with C-index equal to zero are local minimizers of MPCC. For nonvanishing C-indices we obtain all kinds of saddle points for MPCC. Overall, the C-index uniquely determines the topological type of a nondegenerate C-stationary point.

Finally, we define the following auxiliary index subsets, which depend on the signs of multipliers $(\bar \sigma, \bar \rho)$ corresponding to a W-stationary point $\bar x \in M$:
% W- or C- ?????
\[
\begin{array}{ll}
a_{01}^-\left(\bar x\right)=\left\{j \in a_{01}\left(\bar x\right)\,\left\vert\,  \bar \sigma_{1,j} <0 \right.\right\}, & a_{10}^-\left(\bar x\right)=\left\{j \in a_{10}\left(\bar x\right)\,\left\vert\,  \bar \sigma_{2,j} <0 \right.\right\}, \\
a_{01}^0\left(\bar x\right)=\left\{j \in a_{01}\left(\bar x\right)\,\left\vert\,  \bar \sigma_{1,j} =0 \right.\right\}, & a_{10}^0\left(\bar x\right)=\left\{j \in a_{10}\left(\bar x\right)\,\left\vert\,  \bar \sigma_{2,j} =0 \right.\right\},\\
a_{01}^+\left(\bar x\right)=\left\{j \in a_{01}\left(\bar x\right)\,\left\vert\,  \bar \sigma_{1,j} >0 \right.\right\}, & a_{10}^+\left(\bar x\right)=\left\{j \in a_{10}\left(\bar x\right)\,\left\vert\,  \bar \sigma_{2,j} >0 \right.\right\},\\
\end{array}
\]
\[
\begin{array}{l}
a_{00}^-\left(\bar x\right)=\left\{j \in a_{00}\left(\bar x\right)\,\left\vert\,  \bar \varrho_{1,j},\bar \varrho_{2,j} <0 \right.\right\},\\
a_{00}^0\left(\bar x\right)=\left\{j \in a_{00}\left(\bar x\right)\,\left\vert\,  \bar \varrho_{1,j} \cdot \bar \varrho_{2,j} =0 \right.\right\},
\\
a_{00}^+\left(\bar x\right)=\left\{j \in a_{00}\left(\bar x\right)\,\left\vert\,  \bar \varrho_{1,j},\bar \varrho_{2,j} >0 \right.\right\}.
\end{array}
\]
%%%%%%%%%%%%%%%%%%%%%%%%%%%%%%%%%%%%%%
%Scholtes
%%%%%%%%%%%%%%%%%%%%%%%%%%%%%%%%%%%%%%
\subsection{Preliminaries on the Scholtes regularization}\label{sec:scholtes}

 Obviously, the Scholtes regularization $\mathcal{S}$ 
falls into the scope of nonlinear programming (NLP). For the readers' convenience, let us briefly apply the NLP theory, see e.g.~\cite{jongen:2000}, to $\mathcal{S}$. For that, the index sets of active constraints associated with $x \in M^\mathcal{S}$ will be denoted by 
\[
\mathcal{H}(x)=\left\{j\,\left\vert \, F_{1,j}(x)\cdot F_{2,j}(x) = t \right.\right\}, \]
\[
\mathcal{N}_1(x)=\left\{j\,\left\vert \, F_{1,j}(x)= 0 \right.\right\},\quad \mathcal{N}_2(x)=\left\{j\,\left\vert \, F_{2,j}(x)= 0 \right.\right\}.
\]

%%%%%%%%%%%%%%%%%%%%%%%%%%%%%%%%%%%%%%
%KKT
%%%%%%%%%%%%%%%%%%%%%%%%%%%%%%%%%%%%%%

We start with the standard linear independence constraint qualification.

\begin{definition} [LICQ]
We say that a feasible point $x \in M^\mathcal{S}$ of $\mathcal{S}$ satisfies the linear independence constraint qualification (LICQ) if the following vectors are linearly independent:
\[
F_{2,j}(x) D F_{1,j}\left(x\right) + F_{1,j}\left(x\right) D F_{2,j}\left(x\right),\,j \in \mathcal{H}\left(x\right),\]
\[
D F_{1,j}\left(x\right), j \in \mathcal{N}_1\left(x\right),\quad D F_{2,j}\left(x\right), j \in \mathcal{N}_2\left(x\right).
\]    
\end{definition}

It happens that LICQ for $M^\mathcal{S}$ is inherited by MPCC-LICQ for $M$.

\begin{lemma}[MPCC-LICQ vs.~LICQ, \cite{scholtes:2001}]
\label{lem:LICQ-MPCCLICQ}
Let a feasible point $\bar x \in M$ of MPCC fulfill MPCC-LICQ. Then, LICQ holds at all feasible points $x \in M^\mathcal{S}$ of $\mathcal{S}$ for all sufficiently small $t$, whenever they are sufficiently close to $\bar x$. 
\end{lemma}

The topologically relevant stationarity notion for NLP is that of Karush-Kuhn-Tucker points.
%
% Definition with active sets as above???
%
\begin{definition} [Karush-Kuhn-Tucker point]
       A feasible point $x \in M^\mathcal{S}$ of $\mathcal{S}$ is called Kurush-Kuhn-Tucker point if there exist multipliers 
\[
\eta_j,\nu_{1,j}, \nu_{2,j}, \,j\in\left\{1,\ldots,\kappa\right\},
\]
such that the following conditions hold:
\begin{equation}
   \label{eq:kkt-1} 
   \begin{array}{rcl}
\displaystyle
D f\left( x\right) &=& 
-\displaystyle\sum\limits_{j=1}^{\kappa} 
\eta_j \left(F_{2,j}(x)D F_{1,j}\left(x\right) + F_{1,j}\left(x\right) D F_{2,j}\left(x\right)\right)\\&&
+ \displaystyle\sum\limits_{j=1}^{\kappa}   \nu_{1,j} D F_{1,j}\left(x\right)
+ \displaystyle\sum\limits_{j=1}^{\kappa}   \nu_{2,j} D F_{2,j}\left(x\right).
\end{array}
\end{equation}
\begin{equation}
   \label{eq:kkt-2} 
   \begin{array}{c}
 \eta_j\left(t-F_{1,j}(x)\cdot F_{2,j}(x)\right)=0, F_{1,j}(x) \cdot F_{2,j}(x)\le t, \eta_j\ge 0,\\ \\ 
 \nu_{1,j}F_{1,j}(x)=0,F_{1,j}(x) \ge 0, \nu_{1,j}\ge 0, \\ \\ \nu_{2,j}F_{2,j}(x)=0,F_{2,j}(x)\ge 0, \nu_{2,j}\ge 0, j\in\{1,\ldots,\kappa\}.
    \end{array}
\end{equation}
\end{definition}
For a Karush-Kuhn-Tucker point $x \in M^\mathcal{S}$ with multipliers $(\eta, \nu)$ -- which are unique under LICQ -- it is convenient to define the appropriate Lagrange function as
\[
 \begin{array}{rcl}
L^{\mathcal{S}}(x)&=&f(x) - \displaystyle\sum\limits_{j=1}^{\kappa} 
\eta_j \left(t-F_{1,j}(x)F_{2,j}\left(x\right)\right)-
\displaystyle\sum\limits_{j=1}^{\kappa}  \nu_{1,j} F_{1,j}\left(x\right)
- \displaystyle\sum\limits_{j=1}^{\kappa} \nu_{2,j} F_{2,j}\left(x\right).
\end{array}
\]
The corresponding tangent space is given by
\[
    \mathcal{T}^\mathcal{S}_{x}=\left\{
\xi \in \R^{n}\,\left\vert\, \begin{array}{l}
\left(F_{2,j}\left(x\right)DF_{1,j}\left(x\right)+F_{1,j}\left(x\right)DF_{2,j}\left(x\right)\right)\xi=0, j \in \mathcal{H}\left(x\right),\\
DF_{1,j}\left(x\right)\xi=0, j \in \mathcal{N}_1\left(x\right),
DF_{2,j}\left(x\right)\xi=0, j \in \mathcal{N}_2\left(x\right)
\end{array}
\right.\right\}.
\]

Now, we are ready to recall the notion of nondegeneracy for Karush-Kuhn-Tucker points. 

\begin{definition}[Nondegenerate Karush-Kuhn-Tucker point]
A Karush-Kuhn-Tucker point $x \in M^\mathcal{S}$ of $\mathcal{S}$ with multipliers $(\eta, \nu)$
is called nondegenerate if 

ND1: LICQ holds at $x$,

ND2: the multipliers for active inequality constraints do not vanish, i.e.~$\eta_j>0$,\,$j \in \mathcal{H}(x)$, $\nu_{1,j} >0$,\,$j \in \mathcal{N}_1(x)$, and $\nu_{2,j} >0$,\,$j \in \mathcal{N}_2(x)$,

ND3: the restricted Hessian matrix $D^2 L^\mathcal{S}(x)\restriction_{\mathcal{T}^\mathcal{S}_{x}}$ is nonsingular.
\end{definition}

The dimension of the cell to be attached in the framework of Morse theory for NLP corresponds to the quadratic index of a nondegenerate Karush-Kuhn-Tucker point.

\begin{definition} [Quadratic index]
    Let $x \in M^\mathcal{S}$ be a nondegenerate Karush-Kuhn-Tucker point of $\mathcal{S}$ with unique multipliers $(\eta, \nu)$.
The  number of negative eigenvalues of the matrix $D^2 L^\mathcal{S}(x)\restriction_{\mathcal{T}^\mathcal{S}_{x}}$ is called its quadratic index ($QI$).
\end{definition}

%%%%%%%%%%%%%%%%%%%%%%%%%%%%%%%%%%%%%%
%Lemmas for Scholtes
%%%%%%%%%%%%%%%%%%%%%%%%%%%%%%%%%%%%%%

We shall eventually use the following relations between the multipliers of MPCC and of $\mathcal{S}$ as the parameter $t$ tends to zero.

\begin{proposition}[Multipliers, \cite{scholtes:2001}]
\label{prop:com:scholtes}
    Suppose a sequence of Karush-Kuhn-Tucker points $x^{t} \in M^\mathcal{S}$ of $\mathcal{S}$ with multipliers $\left(\eta^t, \nu^t\right)$ converges to $\bar x$ for $t \to 0$. Let MPCC-LICQ be fulfilled at $\bar x \in M$. Then, $\bar x$ is a C-stationary point of MPCC with the unique multipliers $(\bar \sigma, \bar \varrho)$ satisfying:
\begin{itemize}
    \item[a)] for $j \in a_{01}\left(\bar x\right)$, $\bar \sigma_{1,j}=\left\{
\begin{array}{ll}
   \lim\limits_{t \to 0} \nu_{1,j}^{t}  &j \not\in I_0,  \\
    -\lim\limits_{t \to 0} \eta_{j}^{t}F_{2,j}\left(x^t\right)  &j \in I_0,
\end{array}\right.$
\item[b)] for $j \in a_{10}\left(\bar x\right)$, $\bar \sigma_{2,j}=\left\{
\begin{array}{ll}
   \lim\limits_{t \to 0} \nu_{2,j}^{t}  &j \not\in I_0,  \\
    -\lim\limits_{t \to 0} \eta_{j}^{t}F_{1,j}\left(x^t\right)  &j \in I_0,
\end{array}\right.$
\item[c)] for $j \in a_{00}\left(\bar x\right)$, $\bar \varrho_{1,j}=\left\{
\begin{array}{ll}
   \lim\limits_{t \to 0} \nu_{1,j}^{t}  &j \not\in I_0,  \\
    -\lim\limits_{t \to 0} \eta_{j}^{t}F_{2,j}\left(x^t\right)  &j \in I_0,
\end{array}\right.$ and $\bar \varrho_{2,j}=\left\{
\begin{array}{ll}
   \lim\limits_{t \to 0} \nu_{2,j}^{t}  &j \not\in I_0,  \\
    -\lim\limits_{t \to 0} \eta_{j}^{t}F_{1,j}\left(x^t\right)  &j \in I_0,
\end{array}\right.$
\end{itemize}    
where \[
I_0=\left\{j\, \left\vert \, j \in \mathcal{H}\left(x^t\right) \mbox{for infinitely many $t$}\right.\right\}.
\]
\end{proposition}

By means of Proposition \ref{prop:com:scholtes}, it is easy to relate the active index (sub)sets of MPCC and $\mathcal{S}$ as the parameter $t$ is sufficiently close to zero.

\begin{corollary}[Active index sets]
\label{cor:activeindex:scholtes}
Under the assumptions of Proposition \ref{prop:com:scholtes}, we have for the active index subsets of the multipliers $(\bar \sigma, \bar \varrho)$ of $\bar x \in M$ for all sufficiently small $t$:
\begin{itemize}
    \item[a)] $a_{01}^-\left(\bar x\right)\subset \mathcal{H}\left(x^t\right)$ and
 $a_{01}^+\left(\bar x\right)\subset \mathcal{N}_1\left(x^t\right)$,
\item[b)] $a_{10}^-\left(\bar x\right)\subset \mathcal{H}\left(x^t\right)$  and
 $a_{10}^+\left(\bar x\right)\subset \mathcal{N}_2\left(x^t\right)$,
\item[c)] $a_{00}^-\left(\bar x\right)\subset \mathcal{H}\left(x^t\right)$  and
 $a_{00}^+\left(\bar x\right)\subset \mathcal{N}_1\left(x^t\right) \cap \mathcal{N}_2\left(x^t\right)$.
\end{itemize}

\end{corollary}

\section{Main results}
\label{sec:main}

We start by examining the convergence behavior of the Scholtes regularization. This is done up to the topological type of C-stationary points of the corresponding MPCC. For that, we find out what happens to the C-index of a C-stationary point which is the limit of a sequence of nondegenerate Karush-Kuhn-Tucker points of $\mathcal{S}$ as $t \rightarrow 0$. First of all, such limiting C-stationary points may well become degenerate, see Example \ref{ex:ndc2fail}. {\color{black} That nondegenerate stationary points may converge to degenerate limiting points is a very general phenomenon, which even occurs in simple examples of unconstrained optimization problems. We point out that this effect is unrelated to the peculiarities of the Scholtes regularization.}

%%%%%%%%%%%% Anomaly Degeneracy

\begin{example} [Failure of NDC2]
\label{ex:ndc2fail}
We consider the following Scholtes regularization $\mathcal{S}$ with $n=2$ and $\kappa=1$:
%, $F_{1,1}(x)=x_1$, and $F_{2,1}(x)=x_2$:
\[
\begin{array}{rll}
\mathcal{S}:&\min\limits_{x}&\frac{1}{4}\left(x_1^4+x_2^4\right)-\frac{1}{2}\left(x_1^2+x_2^2\right) \\
&\mbox{s.\,t.}& x_1\cdot x_2\le t, \quad x_1\ge 0, x_2\ge 0,
\end{array}
\]
as well as the point $x^t=\left(\sqrt{t},\sqrt{t}\right)$ for $t\in \left(0,\frac{1}{2}\right)$. The latter is a nondegenerate Karush-Kuhn-Tucker point of $\mathcal{S}$. Indeed, it holds
\[
\begin{pmatrix}
\sqrt{t}^3-\sqrt{t}\\
\sqrt{t}^3-\sqrt{t}\\
\end{pmatrix}
= -\eta^t_1
\begin{pmatrix}
    \sqrt{t}\\ \sqrt{t}
\end{pmatrix},
\]
where $\eta^t=1-t>0$.
Obviously, ND1 and ND2 are fulfilled. Let us show that ND3 holds.
We get for the tangent space at $x^t$
\[
    \mathcal{T}^\mathcal{S}_{x^t}=\left\{
\xi \in \R^{2}\,\left\vert\, 
-\xi_1=\xi_2
\right.\right\}.
\]
For the Hessian matrix of the Lagrange function at $x^t$ we have
\[
D^2L^{\mathcal{S}}\left(x^t \right)=
\begin{pmatrix}
  3\left(x_1^t\right)^2-1 & 0\\
  0 & 3\left(x_2^t\right)^2-1
\end{pmatrix}
+\eta^t 
\begin{pmatrix}
    0 & 1\\
    1 & 0
\end{pmatrix}
=\begin{pmatrix}
    3t-1 & 1-t\\
    1-t & 3t-1
\end{pmatrix}.
\]
Consequently, we have for $\xi \in \mathcal{T}^\mathcal{S}_{x^t}$ with $\xi\ne 0$:
\[
\xi^T D^2L^{\mathcal{S}}\left(x^t \right)\xi =
%(3t-1)\xi_1^2+2(1-t)\xi_1 \xi_2 +(3t-1)\xi_2^2=
%(3t-1)\xi_1^2-2(1-t)\xi_1^2 +(3t-1)\xi_1^2=
(8t-4)\xi_1^2\ne 0.
\] 
Further, the limiting point of the sequence $x^t$ for $t\to 0$ is $\bar x=(0,0)$. However, $\bar x$ is a degenerate C-stationary point of the corresponding MPCC. This is due to the violation of NDC2, since it holds:
\[
\begin{pmatrix}
    0\\0
\end{pmatrix}=
\bar \varrho_{1,1}
\begin{pmatrix}
1\\ 0    
\end{pmatrix}
+\bar \varrho_{2,1}
\begin{pmatrix}
    0\\1
\end{pmatrix}
\]
with the unique multipliers $\bar \varrho_{1,1}=\bar \varrho_{2,1}=0$.
\qed
\end{example}

%%%%%%%%%%%%%%%%%%%%%%%%%%%%%%%%%%%%%%
%Theorem KKT--> C-Stat
%%%%%%%%%%%%%%%%%%%%%%%%%%%%%%%%%%%%%%

In view of Example \ref{ex:ndc2fail}, we shall additionally assume that the limiting C-stationary point obtained via Scholtes regularization is nondegenerate. Theorem \ref{thm:kktsequence} {\color{black} and Example \ref{ex:NDC4}} show that it may nevertheless come to a reduction of its C-index if compared with the quadratic index of the approximating Karush-Kuhn-Tucker points{{\color{black}, see Appendix for the corresponding proof}. %In other words, the convergence behavior of the Scholtes regularization is anomalous in the sense that the topological type of the limiting C-stationary point of MPCC may change.

\begin{theorem}[Convergence]
\label{thm:kktsequence}
%Let Assumption \ref{ass:c} be fulfilled and $\varepsilon \le \frac{1}{n-s}$. 
Suppose a sequence of nondegenerate Karush-Kuhn-Tucker points $x^{t} \in M^\mathcal{S}$ of 
$\mathcal{S}$ with quadratic index $q$ converges to $\bar x$ for $t \to 0$. If $\bar x \in M$ is a nondegenerate C-stationary point of MPCC with multipliers $(\bar \sigma, \bar \varrho)$, then we have for its C-index: 
\begin{equation}
    \label{eq:main-index}
\max\left\{q - \left\vert\left\{j\in a_{01}\left(\bar x\right)\,\left\vert\,\bar \sigma_{1,j}=0\right.\right\}\right\vert - \left\vert\left\{j\in a_{10}\left(\bar x\right)\,\left\vert\,\bar \sigma_{2,j}=0\right.\right\} \right\vert, 0\right\}\le CI \leq q.
\end{equation}
If additionally NDC4 holds at $\bar x$, then the indices coincide, i.e.~$CI= q$.
\end{theorem}

%%%%%%%%%%%%%%%%%%%%%%%%%%%%%%%%%%%%%%
%Example Index decrease
%%%%%%%%%%%%%%%%%%%%%%%%%%%%%%%%%%%%%%

We emphasize that in absence of NDC4 the shift of the C-index in Theorem \ref{thm:kktsequence} cannot be avoided in general. This is illustrated by Example \ref{ex:NDC4}, where the lower bound for $CI$ in (\ref{eq:main-index}) is attained.  

\begin{example}[Necessity of NDC4]
\label{ex:NDC4}
We consider the following Scholtes regularization $\mathcal{S}$ with
$n=3$ and $\kappa=3$: 
%and 
%$F_{1,1}(x)=x_1, F_{2,1}(x)=x_3-2x_2, F_{1,2}(x)=x_1+x_2, F_{2,2}(x)=2-x_3, F_{1,3}(x)=x_3-1, F_{2,3}(x)=x_1-x_2+x_3$:
\[
\begin{array}{rll}
 \mathcal{S}:& \min\limits_x & -x_1-x_2+2x_1x_2+x_3^2\\
 &\text{s.t.}& x_1\cdot (x_3-2x_2) \le t, x_1\ge 0, x_3-2x_2 \ge 0,\\
 &&(x_1+x_2)\cdot (2-x_3) \le t, x_1+x_2\ge 0, 2-x_3\ge 0,\\
 &&(x_3-1)\cdot (x_1-x_2+x_3) \le t, x_3-1 \ge 0, x_1-x_2+x_3 \ge 0.
\end{array}
\]
Let us show that $x^t=\left(\frac{t}{2}, \frac{t}{2},1\right)$ is a nondegenerate Karush-Kuhn-Tucker point of $\mathcal{S}$ if $0<t<1$. Indeed, it is feasible with $\mathcal{H}\left(x^t\right)=\{2\}$, $\mathcal{N}_1\left(x^t\right)=\{3\}$, $\mathcal{N}_2\left(x^t\right)=\emptyset$, and it holds
\[
\begin{pmatrix}
    -1+t\\-1+t\\2
\end{pmatrix}
=
-\eta^t_{2}\begin{pmatrix}
    1\\1\\-t
\end{pmatrix}
+\nu^t_{1,3}\begin{pmatrix}
    0\\0\\1
\end{pmatrix}
\]
with the unique multipliers $\eta^t_{2}=1-t, \nu^t_{1,3}=2-t+t^2$.
Obviously, ND1 and ND2 are fulfilled.
Further, for the tangent space at $x^t$ we have
\[
\mathcal{T}^{\mathcal{S}}_{x^t}=\left\{\xi \in \R^3\,\left\vert \,\xi_1=-\xi_2,\xi_3=0\right.\right\}.
\]
For the Hessian matrix of the Lagrange function at $x^t$ we get
\[
D^2L^{\mathcal{S}}\left(x^t\right)=
\begin{pmatrix}
    0&2&t-1\\
    2&0&t-1\\
    t-1&t-1&2
\end{pmatrix}.
\]
Thus, for $\xi \in \mathcal{T}^{\mathcal{S}}_{x^t}$ with $\xi \ne 0$ it holds
\[
\xi^T D^2L^{\mathcal{S}}\left(x^t\right) \xi = 
4\xi_1 \xi_2 +2(t-1) \xi_1 \xi_3 +2(t-1)\xi_2 \xi_3 +2\xi_3^2=-4\xi_2^2<0.
\]
We conclude that ND3 is also fulfilled and, thus, $x^t$ is nondegenerate. For the quadratic index of $x^t$ we have $q=1$, since the matrix $D^2 L^\mathcal{S}(x^t) \restriction_{\mathcal{T}^\mathcal{S}_{x^t}}$ is negative definite and $\mbox{dim}\left(\mathcal{T}^{\mathcal{S}}_{x^t}\right)=1$.

Further, the sequence $x^t$ converges to $\bar x=(0,0,1)$ for $t \to 0$. Therefore, let us focus on $\bar x$ as well as the underlying 
\[
\begin{array}{rll}
 \mbox{MPCC:}& \min\limits_x & -x_1-x_2+2x_1x_2+x_3^2\\
 &\text{s.t.}& x_1\cdot (x_3-2x_2) =0, x_1\ge 0, x_3-2x_2 \ge 0,\\
 &&(x_1+x_2)\cdot (2-x_3) =0, x_1+x_2\ge 0, 2-x_3\ge 0,\\
 &&(x_3-1)\cdot (x_1-x_2+x_3) =0, x_3-1 \ge 0, x_1-x_2+x_3 \ge 0.
\end{array}
\]
In view of Proposition \ref{prop:com:scholtes}, $\bar x$ is a C-stationary point. 
%index sets $a_{01}\left(\bar x\right)=\{1,2,3\}, a_{10}\left(\bar x\right)=\emptyset$, and $a_{00}\left(\bar x\right)=\emptyset$ as well as
It holds
\[
\begin{pmatrix}
    -1\\-1\\2
\end{pmatrix}
=
\bar \sigma_{1,1}\begin{pmatrix}
    1\\0\\0
\end{pmatrix}
+\bar \sigma_{1,2}\begin{pmatrix}
    1\\1\\0
\end{pmatrix}
+\bar \sigma_{1,3}\begin{pmatrix}
   0\\0\\1
\end{pmatrix}
\] with the unique multipliers $\bar \sigma_{1,1}=0$, $\bar \sigma_{1,2}=-1$, and $\bar \sigma_{1,3}=2$.
It easily follows that $\bar x$ satisfies conditions NDC1, NDC2, and NDC3. Hence, $\bar x$ is nondegenerate. For the quadratic and biactive index of $\bar x$ we have $QI=0$ and $BI=0$, respectively.
For the C-index of $\bar x$ it holds then $CI=QI+BI=0$.
Overall, the quadratic index of $x^t$ is strictly greater than the C-index of $\bar x$ and the lower bound from Theorem \ref{thm:kktsequence} is attained, i.e.~
\[
\max\left\{q - \left\vert\left\{j\in a_{01}\left(\bar x\right)\,\left\vert\,\bar \sigma_{1,j}=0\right.\right\}\right\vert - \left\vert\left\{j\in a_{10}\left(\bar x\right)\,\left\vert\,\bar \sigma_{2,j}=0\right.\right\} \right\vert, 0\right\}=0= CI<q=1.
\]
Here, the sequence $x^t$ of saddle points of $\mathcal{S}$ converges to the minimizer $\bar x$ of MPCC. This is caused by to the violation of NDC4 at $\bar x$.
\qed
\end{example}

%%%%%%%%%%%%%%%%%%%%%%%%%%%%%%%%%%%%%%
%Theorem Existens KKT
%%%%%%%%%%%%%%%%%%%%%%%%%%%%%%%%%%%%%%

Next, we turn our attention to the well-posedness of Scholtes relaxation in the vicinity of C-stationary points of the corresponding MPCC{\color{black}, see Appendix for the corresponding proof}. Also here the validity of NDC4 becomes crucial.

\begin{theorem}[Well-posedness]
\label{thm:wellposedness}
   Let $\bar x \in M$ be a nondegenerate C-stationary point of MPCC with C-index $c$, additionally, fulfilling NDC4. Then, for all sufficiently small $t$ there exists a nondegenerate Karush-Kuhn-Tucker point $x^t \in M^\mathcal{S}$ of $\mathcal{S}$ within a neighborhood of $\bar x$, which has the same quadratic index $c$.
   Moreover, for any fixed $t$ sufficiently small, such $x^t$ is the unique Karush-Kuhn-Tucker point of $\mathcal{S}$ in a sufficiently small neighborhood of $\bar x$.
\end{theorem}

We {\color{black} proceed} with a genericity result for NDC4. 

\begin{proposition}[NDC4 is generic]
\label{prop:generic}
 The subset of MPCC defining functions $f, F_1, F_2$, for which each nondegenerate
C-stationary point of a corresponding MPCC fulfills NDC4, is $C^2_s$-open and -dense.
%    
%    Let $\mathcal{F} \subset C^2\left(\R^n,\R\right) \times C^2\left(\R^n,\R^{\kappa}\right) \times C^2\left(\R^n,\R^{\kappa}\right)$ be the subset of MPCC defining functions for which each nondegenerate C-stationary point fulfills NDC4. Then, $\mathcal{F}$ is $C^2_s$-open and -dense.
\end{proposition}

\proof    Let us fix index sets $A_1, A_2 \subset \left\{1,\ldots,\kappa\right\}$ of
active defining functions $F_1$ and $F_2$, respectively. Furthermore, let us fix subsets of those:
$B_1\subset A_1$ and $B_2\subset A_2$.
For this choice we consider the set  $M_{A_1,A_2,B_1,B_2}$ of $x \in \R^n$ such that the following conditions are satisfied:
\begin{itemize}
    \item[](c1)\; $F_{1,j}=0$ for all $j \in A_1$, and $F_{1,j}\ne 0$ for all $j \not\in A_1$,
    \item[](c2)\; $F_{2,j}=0$ for all $j \in A_2$, and $F_{2,j}\ne 0$ for all $j \not\in A_2$,
    \item[](c3)\; $D f(x) \in \mbox{span}\left\{D F_{1,j}(x), j \in A_1\backslash B_1,
    D F_{2,j}(x), j \in A_2\backslash B_2\right\}$.
\end{itemize}
Here, (c1) and (c2) refer to active constraints, while (c3) mimics C-stationarity and possible violations of NDC2 and NDC4.
We show that $M_{A_1,A_2,B_1,B_2}$ is generically empty whenever one of the sets $B_1$ or $B_2$ is nonempty.
For that, let us consider the codimension of  $M_{A_1,A_2,B_1,B_2}$. We observe, that the available degrees of freedom of the variables involved are $n$. Conditions (c1) and (c2) cause a loss of freedom of $|A_1|$, and
$|A_2|$, respectively.
Further (c3) causes an additional loss of freedom of $n-\left|A_1\right|+\left|B_1\right|-\left|A_2\right|+\left|B_2\right|$ due to NDC1. Hence, the total loss of freedom is $n+\left|B_1\right|+\left|B_2\right|$.
We conclude that a violation of NDC2 or NDC4, i.e.~$B_1 \not = \emptyset$ or $B_2 \not = \emptyset$, would imply that the total available degrees of freedom $n$ are exceeded. Generically, the sets $M_{A_1,A_2,B_1,B_2}$ must be thus empty by virtue of the jet
transversality theorem from \cite{jongen:2000}.

The openness result follows by a standard reasoning, using the fact that, locally C-stationarity can be written
via stable equations. Then, the implicit function theorem for Banach spaces can be applied to
follow C-stationary points with respect to (local) $C^2$-perturbations of defining functions. Finally,
a standard globalization procedure exploiting the specific properties of the strong $C^2_s$-topology can be used to construct a (global) $C^2_s$-neighborhood of problem data for which the nondegeneracy
property is stable, cf.~\cite{jongen:2000}. \qed

{\color{black} Finally, we comment on NDC4 from the topological perspective, i.e. in view of the typical cell-attachment results from the Morse theory for MPCC.}

{\color{black}
\begin{remark}[NDC4 from the topological point of view]
As it became clear from Theorems \ref{thm:kktsequence} and \ref{thm:wellposedness}, NDC4 is essential for guaranteeing stability of the Scholtes regularization if focusing on the topological type of the corresponding stationary points.
However, we would like to stress that NDC4 is not of topological significance for MPCC. Indeed, the non-biactive multipliers play a role neither in the original definition of nondegeneracy for C-stationary points, i.e. NDC1-NDC3, nor in the definition of their C-index, cf. Definition \ref{def:c-index}. 
In particular, negative signs of non-biactive multipliers do not contribute to the C-index. More than that, such an adjustment would be just inappropriate from the topological point of view. This is since the C-index as defined in \cite{ralph:2011} has been shown to correspond to the dimension of a cell to be attached in order to describe the topological changes of the MPCC lower level sets within the scope of the Morse theory, see \cite{jongen:2009}. Additionally, we emphasize that non-biactive multipliers are taken with respect to the local equality constraints in MPCC and, as it is well-known from the Morse theory for NLP \cite{jongen:2000}, the multiplier information for equality constraints is therefore irrelevant.
\qed
\end{remark}
}

%Note that nondegeneracy is a generic property, see \cite{jongen:2009}. 
%Hence, the assumption of nondegeneracy could also be omitted in Proposition \ref{prop:generic}.
    
\section{Comparison to the existing literature}
\label{sec:compare}

The aim of this section is to relate our results in Theorems \ref{thm:kktsequence} and \ref{thm:wellposedness} to those known from the literature. We refer mainly to the references \cite{scholtes:2001, still:2007, ralph:2004}.

\subsection{Results of \citet{scholtes:2001}}

%%%%%% Hinrichtung: M-stationär + SONC, Zusammenhang zu Theorem 1: C_x vs. T_x, and SONC vs. NDC3 under NDC2/NDC4, keine Sattelpunkte

%%%%%% Rückrichtung + NDC4 über strenge Stabilität mit Eindeutigkeit, Spezialfall von Theorem 2, aber nur für Index 0 

In order to adequately state the results from \cite{scholtes:2001},
we need to introduce some auxiliary notions for the Scholtes regularization.
For a Karush-Kuhn-Tucker point $x \in M^\mathcal{S}$ of $\mathcal{S}$ with the corresponding multipliers $(\eta, \nu)$ we define the following index sets:
\[
\mathcal{H}^+\left(x\right)=\left\{
j \in \mathcal{H}\left(x\right) 
\,\left\vert \, \eta_{j}>0
\right. \right\},
\]
\[
\mathcal{N}^+_1\left(x\right)=\left\{
j \in \mathcal{N}_1\left(x\right) 
\,\left\vert \, \nu_{1,j}>0
\right. \right\}, \quad
\mathcal{N}^+_2\left(x\right)=\left\{
j \in \mathcal{N}_2\left(x\right)
\,\left\vert \, \nu_{2,j}>0
\right. \right\}.
\]
Additionally, we define the corresponding set of critical directions:
    \[
\mathcal{C}^{\mathcal{S}}_{x}=\left\{\xi \in \R^n \,\left\vert \,
\begin{array}{l}
\left(F_{2,j}\left(x\right)DF_{1,j}\left(x\right)+F_{1,j}\left(x\right)DF_{2,j}\left(x\right)\right)\xi=0, j \in \mathcal{H}^+\left(x\right), \\
DF_{1,j}\left(\bar x\right)\xi=0, j \in \mathcal{N}^+_1\left(x\right), \\
DF_{2,j}\left(\bar x\right)\xi=0, j \in \mathcal{N}^+_2\left(x\right)
\end{array}
\right.\right\}. 
\]
Next, we recall a second order necessary condition for the Scholtes relaxation.

\begin{definition}[SONC, \cite{scholtes:2001}, \cite{fletcher:1987}]
    A Karush-Kuhn-Tucker point $x$ of $\mathcal{S}$ with multipliers $\left(\eta, \nu\right)$ is said to fulfill  
    %with the unique multipliers $(\bar \sigma, \bar \varrho)$. We say that it fulfills 
    the second order necessary condition (SONC) if for every $\xi \in {\mathcal{C}}^{\mathcal{S}}_{x}$ it holds: 
    \[
    \xi^T D^2L^{\mathcal{S}}\left(
    x\right)\xi \ge 0.
    \]
\end{definition}

\begin{comment}
\begin{definition}[ULSC, \cite{scholtes:2001}]
    A W-stationary point $\bar x$ of MPCC is said to fulfill  
    %with the unique multipliers $(\bar \sigma, \bar \varrho)$. We say that it fulfills 
    the upper level strict complementarity (ULSC) if there exist multipliers $(\bar \sigma, \bar \varrho)$ with
    \[
    \varrho_{1,j} \cdot \varrho_{2,j}\ne 0 \mbox{ for all } j \in a_{00}\left(\bar x\right).
    \]
\end{definition}  
\end{comment}
%Note that for a C-stationary point $\bar x$ of MPCC, NDC2 and ULSC are equivalent.
%However, if $\bar x$ is also C-stationary, then clearly NDC2 and ULSC coincide. 

Now, we are ready to present what has been altogether shown in \cite{scholtes:2001} on the convergence properties of the Scholtes relaxation, see Theorems 3.1, 3.3, and Corollary 3.4 there.

\begin{theorem}[Convergence, \cite{scholtes:2001}]
\label{thm:scholtes-1}
  Suppose a sequence of Karush-Kuhn-Tucker points $x^{t} \in M^\mathcal{S}$ of $\mathcal{S}$ with multipliers $\left(\eta^t, \nu^t\right)$ converges to $\bar x$ for $t \to 0$. Assume SONC holds at each $x^t$ and let MPCC-LICQ be fulfilled at $\bar x \in M$. Then, $\bar x$ is M-stationary.
  If additionally $\bar x$ fulfills NDC2, then $\bar x$ is S-stationary.
\end{theorem}

    Let us compare Theorems \ref{thm:kktsequence} and \ref{thm:scholtes-1}. 
First, it is clear that Theorems \ref{thm:kktsequence} and \ref{thm:scholtes-1} are of independent interest. Namely, in Theorem \ref{thm:kktsequence} we can handle the cases where SONC is violated, but ND3 is fulfilled. This includes nondegenerate saddle points of $\mathcal{S}$, which are characterized by the nonvanishing quadratic index. On the other side, Theorem \ref{thm:scholtes-1} does not need ND2, i.e.~the strict complementarity for $\mathcal{S}$. This can be explained by the fact that SONC is stated in terms of critical rather than tangential directions, cf.~ND3. Actually, it is more important that the assumption of NDC2 for MPCC can be avoided in Theorem \ref{thm:scholtes-1}. Even without NDC2, M-stationary points are obtained in the limit. Note that they form a proper subset of C-stationary points, which are under consideration in our Theorem \ref{thm:kktsequence}. Again, this stronger result of Theorem \ref{thm:scholtes-1} is due to the validity of SONC, which prevents from considering the Karush-Kuhn-Tucker points with nonvanishing quadratic index. 
Second, Theorems \ref{thm:kktsequence} and \ref{thm:scholtes-1} have similarities if focusing on the case of minimizers. With the additional assumption of NDC2 even S-stationary points can be obtained in the limit, see Theorem \ref{thm:scholtes-1}. In our terms, this corresponds to the vanishing biactive index. Theorem \ref{thm:kktsequence} traces both, the biactive and the quadratic parts of the C-index. If starting with nondegenerate local minimizers of $\mathcal{S}$, the C-index remains zero, see the inequality in Theorem \ref{thm:kktsequence}. Hence, we get local minimizers of MPCC in the limit. Overall, we conclude that our results complement those of \cite{scholtes:2001}. In particular, we discuss the saddle points of $\mathcal{S}$ and, moreover, derive an inequality for the C-index shift while taking the limit in the Scholtes regularization. The latter phenomenon is new as highlighted in Example \ref{ex:NDC4}.

Now, we turn our attention to the well-posedness results of \cite{scholtes:2001}. For that, let us consider an S-stationary point $\bar x \in M$ with the corresponding multipliers $(\bar \sigma, \bar \varrho)$. For the latter we define the following index subsets corresponding to the biactive multipliers:
\[
a_{00}^1(\bar x)=\left\{
j \in a_{00}\left(\bar x\right) 
\,\left\vert \, \bar \varrho_{1,j}>0
\right. \right\}, \quad
a_{00}^2(\bar x)=\left\{
j \in a_{00}\left(\bar x\right) 
\,\left\vert \, \bar \varrho_{2,j}>0
\right. \right\}.
\]
Analogously to above, we define the set of critical directions:
    \[
\mathcal{C}_{\bar x}=\left\{\xi \in \R^n \,\left\vert \,
\begin{array}{l}
DF_{1,j}\left(\bar x\right)\xi=0, j \in a_{01}^-(\bar x) \cup a_{01}^+(\bar x)\cup a_{00}^1(\bar x),\\
DF_{2,j}\left(\bar x\right)\xi=0, j \in a_{10}^-(\bar x) \cup a_{10}^+(\bar x)\cup a_{00}^2(\bar x)
\end{array}
\right.\right\}. 
\]
Also, a tailored second order necessary condition for MPCC can be formulated.

\begin{definition}[SSOSC, \cite{scholtes:2001}]
    Let MPCC-LICQ hold at an S-stationary point $\bar x \in M$ of MPCC with the unique multipliers $(\bar \sigma, \bar \varrho)$. We say that it fulfills the strong second order sufficiency condition (SSOSC) if for every $\xi \in {\mathcal{C}}_{\bar x} \backslash \{0\}$ it holds: 
    \[
    \xi^T D^2L\left(\bar x\right)\xi >0.
    \]
\end{definition}

Let us now cite Theorem 4.1 from \cite{scholtes:2001}, which we state for convenience in terms of our nondegeneracy conditions.

\begin{theorem}[Well-posedness, \cite{scholtes:2001}]
\label{thm:scholtes-2}
    Let $\bar x \in M$ be an S-stationary point of MPCC fulfilling MPCC-LICQ, NDC2, SSOSC, and NDC4.
    Then, for all sufficiently small $t$ there exists a Karush-Kuhn-Tucker point $x^t \in M^\mathcal{S}$ of $\mathcal{S}$ within a neighborhood of $\bar x$.
   Moreover, for any fixed $t$ sufficiently small, such $x^t$ is the unique Karush-Kuhn-Tucker point of $\mathcal{S}$ in a sufficiently small neighborhood of $\bar x$. Further, $x^t$ fulfills the second order sufficient condition.
   %Then, there exists a neighborhood $U$ of $\bar x$, a $\bar t>0$, and a piecewise smooth function $z:[-\bar t, \bar t]\to U$ such that $x(t)=z(t)$ is the unique Karush-Kuhn-Tucker point $x^t$ of $\mathcal{S}(t)$ in $U$ for every $t\in [o, \bar t)$.
\end{theorem}

Note that Theorem \ref{thm:wellposedness} can be seen as a generalization of Theorem \ref{thm:scholtes-2}.
Instead of $\bar x$ being C-stationarity and fulfilling NDC3, it assumes S-stationarity and SSOSC, respectively.
However, C-stationarity is implied by S-stationarity, while NDC3 is implied by SSOSC. The latter is in view of the fact that under NDC2 and NDC4 the tangent space $\mathcal{T}_{\bar x}$ and the set of critical directions $\mathcal{C}_{\bar x}$ coincide. More precisely, under the assumptions of Theorem \ref{thm:scholtes-2}, the biactive index of $\bar x$ vanishes due to S-stationarity. Further, the quadratic index of $\bar x$ must be zero in view of SSOSC. We conclude that Theorem \ref{thm:scholtes-2} coincides with Theorem \ref{thm:wellposedness} in case of a vanishing C-index, thus, if considering local minimizers. Hence, the contribution of Theorem \ref{thm:wellposedness} is that we expand the assertion of Theorem \ref{thm:scholtes-2} to the case of saddle points as well. Another novelty is that we showed the well-posedness of Scholtes regularization up to the topological type of the corresponding stationary points. Not only do we know that in the vicinity of a C-stationary point of MPCC there is a unique Karush-Kuhn-Tucker point of $\mathcal{S}$, but also that the quadratic index of the latter corresponds to the C-index of the former.

\subsection{Results of \citet{ralph:2004}}

%%%%%% Bedingungen: NDC1-NDC4 für Minima impliziert Ralph, aber NDC1-NDC3 und Ralph sind unabhängig; Scholtes impliziert Ralph (echte Verallgemeinerung); Ralph impliziert strenge Stabilität der S-stationären Punkte, aber nicht umgekehrt (Determinantenbedingung). 

%%%%%% Rückrichtung: Aus Ralph folgt Existenz und Eindeutigkeit + strenge Stabilität (?), aber degenerierte KKT-Punkte möglich (Gegenbeispiel: NDC1-NDC3, aber ohne NDC4).

%%%%%%%%%%%%%%%%%%%%%%%%%

% Satz Hinrichtung: Welche Invarianten, z.B. Metric signature or Codimension of a singularity?

% Satz Rückrichtung: In der Nähe von einem streng stabilen C-stationären Punkt (plus noch etwas, NDC4 ?) gibt es (lokal) eindeutige streng stabile KKT-Punkte. (Kontinuumsgegenbeispiel ?)

%%%%%%%%%%%%%%%%%%%%%%%%%%%%%%%%%%%%%%
%Example Degenerate under MPCC-SSOSC
%%%%%%%%%%%%%%%%%%%%%%%%%%%%%%%%%%%%%%

We proceed by a comparison of our results in Section \ref{sec:main} with those of \cite{ralph:2004}.
For that, let us cite Theorem 3.7 of the latter.

\begin{theorem}[Well-posedness, \cite{ralph:2004}]
\label{thm:ralph}
    Let $\bar x \in M$ be an S-stationary point of MPCC fulfilling MPCC-LICQ and SSOSC.
    Then, for all sufficiently small $t$ there exists a Karush-Kuhn-Tucker point $x^t \in M^\mathcal{S}$ of $\mathcal{S}$ within a neighborhood of $\bar x$.
   Moreover, for any fixed $t$ sufficiently small, such $x^t$ is the unique Karush-Kuhn-Tucker point of $\mathcal{S}$ in a sufficiently small neighborhood of $\bar x$.
   %Then, there exists a neighborhood $U$ of $\bar x$, a $\bar t>0$, and a piecewise smooth function $z:[-\bar t, \bar t]\to U$ such that $x(t)=z(t)$ is the unique Karush-Kuhn-Tucker point $x^t$ of $\mathcal{S}(t)$ in $U$ for every $t\in [o, \bar t)$.
\end{theorem}

It is obvious that Theorem \ref{thm:ralph} is a generalization of Theorem \ref{thm:scholtes-2}. Namely, it drops the assumptions of NDC2 and NDC4, while the assertion remains essentially unchanged. 
Unlike Theorem \ref{thm:scholtes-2}, it is however not a special case of our Theorem \ref{thm:wellposedness} anymore.
The reason is that the weaker assumptions of Theorem \ref{thm:ralph} do not guarantee that the unique Karush-Kuhn-Tucker point of $\mathcal{S}$ is nondegenerate. For example, the strict complementarity, i.e.~ND2, might be violated, cf.~Example \ref{ex:ssosc}. 
%Once again, an anomaly of the Scholtes regularization for MPCC can be observed if focusing on the nondegeneracy conditions. 

%%%%%%%%%%%% Anomaly Degeneracy

    \begin{example}[Degeneracy under SSOSC]
\label{ex:ssosc}
    We consider the following MPCC with
$n=2$ and $\kappa=1$: 
\[
\begin{array}{ll}
\min\limits_x & x_1^2+x_2^2\\
\text{s.t.}& x_1\cdot (1-x_2) = 0, x_1\ge 0, 1-x_2 \ge 0.
\end{array}
\]
We show that $\bar x=(0,0)$ is an S-stationary point. Indeed, it holds
\[
\begin{pmatrix}
    0\\0
\end{pmatrix}
=
\bar \sigma_{1,1}\begin{pmatrix}
    1\\0
\end{pmatrix}
\]
with the unique multiplier $\bar \sigma_{1,1}=0$.
It is straightforward to verify that it is nondegenerate, i.e. NDC1-NDC3 hold, and fulfills SSOSC. 
%Since we have for its C-index $c=0$ it is also a minimizer of the considered MPCC. 
However, NDC4 is violated.
Let us focus on the corresponding Scholtes regularization 
\[
\begin{array}{rll}
\mathcal{S}:&\min\limits_x & x_1^2+x_2^2\\
&\text{s.t.}& x_1\cdot (1-x_2) \le t, x_1\ge 0, 1-x_2 \ge 0.
\end{array}
\]
%Note that $\bar x$ fulfills NDC2 and, therefore, PSC. Hence, w
We apply Theorem \ref{thm:ralph} to guarantee the existence of the (locally) unique Karush-Kuhn-Tucker point of $\mathcal{S}$.
Indeed, $x^t=(0,0)$ is this Karush-Kuhn-Tucker point of $\mathcal{S}$ with the
unique multiplier $\nu_{1,1}^t=0$.
However, it violates ND2 and, therefore, $x^t$ is degenerate for $\mathcal{S}$.
\qed
\end{example}

\subsection{Results of \citet{still:2007}}

%%%%%% Rückrichtung in Satz und Remark mit NDC1-NDC3, aber ohne NDC4 über Koordinatentrafo + IFT, Gegenbeispiel für Eindeutigkeit (falsch) + Gegenbeispiel für Index, Gegenbeispiel für ND (Kontinuum).

%%%%%% Frage für später: Ist bei der Rückrichtung NDC1-NDC3 ausreichend für Existenz und/oder Index? Vermutung: NDC1, NDC3 liefern zusammen Existenz; NDC1-NDC3 + in der Nähe ND KKT-Punkte, dann gibt es welche mit dem selben Index.

{\color{black}
In \cite{still:2007}, the authors propose a closely related smoothing approach to MPCC:
\[
\mathcal{S}^=: \quad
\min_{x} \,\, f(x)\quad \mbox{s.\,t.} \quad x \in M^{\mathcal{S}^=}
\]
with
\[
    M^{\mathcal{S}^=}=\left\{x \in\R^n\, \left\vert\,\begin{array}{l}
    F_{1,j}(x) \cdot F_{2,j}(x)= t, F_{1,j}(x) \ge 0,F_{2,j}(x)\ge 0, j=1,\ldots,\kappa
    \end{array} \right.\right\}.
\]
Here, the inequality constraints $F_{1,j}(x) \cdot F_{2,j}(x)\le t, j=1,\ldots,\kappa$, from \cite{scholtes:2001} are replaced by equality constraints.}
%
%Obviously, if $x$ is a Karush-Kuhn-Tucker point (minimizer) of $\mathcal{S}$ with
%$\mathcal{H}\left( x\right)=\left\{1,\ldots,\kappa\right\}$, then it is a Karush-Kuhn-Tucker point (minimizer) of $\mathcal{S}^=$.#
%We note that if $x \in M^{\mathcal{S}}$ is a nondegenerate Karush-Kuhn-Tucker point of $\mathcal{S}$ fulfilling $\mathcal{H}(x)=\{1,\ldots, \kappa\}$, then it holds $x \in M^{\mathcal{S}^=}$. Further, $x$ is a nondegenerate Karush-Kuhn-Tucker point of $\mathcal{S}^=$.
%
In order to state the well-posedness result given in 
\cite{still:2007}, we need some further notions.
%
%\begin{definition}[MPCC-SC, %\cite{still:2007}]
%    An S-stationary point $\bar x \in M$ of MPCC is said to fulfill the strict complementarity slackness for MPCC (MPCC-SC) if it holds:    \[\varrho_{1,j} >0, \varrho_{2,j}>0 \mbox{ for all } j \in a_{00}\left(\bar x\right).\]
%\end{definition}
%In our terms, MPCC-SC corresponds to NDC2 for S-stationary points. 
For an S-stationary points $\bar x \in M$, we define the set of critical directions:
        \[
\mathcal{C}^=_{\bar x}=\left\{\xi \in \R^n \,\left\vert \,
\begin{array}{l}
DF_{1,j}\left(\bar x\right)\xi=0, j \in a_{01}(\bar x) \cup a_{00}^1(\bar x),\\
DF_{2,j}\left(\bar x\right)\xi=0, j \in a_{10}(\bar x) \cup a_{00}^2(\bar x),\\
DF_{1,j}\left(\bar x\right)\xi\ge 0, \mbox{ if } \varrho_{1,j}=0,j \in a_{00}(\bar x),\\
DF_{2,j}\left(\bar x\right)\xi\ge 0, \mbox{ if }  \varrho_{2,j}=0,j \in a_{00}(\bar x)
\end{array}
\right.\right\}. 
\]

\begin{definition}
    An S-stationary point $\bar x \in M$ of MPCC is said to fulfill the second order condition for MPCC (MPCC-SOC) if for every $\xi \in \mathcal{C}^=_{\bar x} \backslash \{0\}$ it holds:
    \[    \xi^T D^2L\left(\bar x\right)\xi >0.\]

\end{definition}
Note that MPCC-SOC is weaker than SSOSC due to $\mathcal{C}^=_{\bar x}\subset{\mathcal{C}}_{\bar x}$.
Having said that, MPCC-SOC is stronger than NDC3 since it holds $\mathcal{T}_{\bar x}\subset \mathcal{C}^=_{\bar x}$.
However, it is also straightforward to verify that for a nondegenerate S-stationary point with vanishing C-index the
conditions NDC3 and MPCC-SOC coincide. If additionally NDC4 holds, then NDC3, MPCC-SOC, and SSOSC coincide respectively. 

%Note that an S-stationary point fulfills MPCC-SC and MPCC-SOC if and only if it fulfills NDC2 and SSOSC.

%Vergleich Still vs. Scholtes

Now, we are ready to cite Theorem 5.1 from \cite{still:2007}. In comparison to Theorem \ref{thm:scholtes-2}, it drops the assumption of NDC4 and only requires a weaker second order condition, namely MPCC-SOC instead of SSOSC, to hold. 

\begin{theorem}[Well-posedness, \cite{still:2007}]
\label{thm:still}
    Let $\bar x \in M$ be a local minimizer of MPCC fulfilling MPCC-LICQ, NDC2, and MPCC-SOC.
    Then, for all sufficiently small $t$ the local minimizers $x^t$ of $\mathcal{S}^=$ (near $\bar x$) are
    uniquely determined.
   %Then, there exists a neighborhood $U$ of $\bar x$, a $\bar t>0$, and a piecewise smooth function $z:[-\bar t, \bar t]\to U$ such that $x(t)=z(t)$ is the unique Karush-Kuhn-Tucker point $x^t$ of $\mathcal{S}(t)$ in $U$ for every $t\in [o, \bar t)$.
\end{theorem}

We recall from \cite{scheel:2000} that any local minimizer of MPCC fulfilling MPCC-LICQ must be S-stationary. It is easy to verify that an S-stationary point fulfills MPCC-LICQ, NDC2, and MPCC-SOC if and
only if it is nondegenerate as a C-stationary point with C-index zero. Thus, the assumptions of Theorem \ref{thm:still} are less strict than those of Theorem \ref{thm:kktsequence}. In particular, NDC4 does not necessarily have to hold. 
\color{black}
The latter is crucial in the context of Scholtes regularization, in order to ensure uniqueness, see Examples \ref{ex:2min}-\ref{ex:continuum2}.

\color{black}
\begin{example}[Non-uniqueness in absence of NDC4]
\label{ex:2min}
    We consider the following Scholtes regularization $\mathcal{S}$ with
$n=3$ and $\kappa=3$: 
%and 
%$F_{1,1}(x)=x_1, F_{2,1}(x)=x_3-2x_2, F_{1,2}(x)=x_1+x_2, F_{2,2}(x)=2-x_3, F_{1,3}(x)=x_3-1, F_{2,3}(x)=x_1-x_2+x_3$:
\[
\begin{array}{rll}
 \mathcal{S}:& \min\limits_x & -x_1-x_2+2x_1x_2+x_3^2\\
 &\text{s.t.}& x_1\cdot (x_3-2x_2) \le t, x_1\ge 0, x_3-2x_2 \ge 0,\\
 &&(x_1+x_2)\cdot (2-x_3) \le t, x_1+x_2\ge 0, 2-x_3\ge 0,\\
 &&(x_3-1)\cdot (x_1-x_2+x_3) \le t, x_3-1 \ge 0, x_1-x_2+x_3 \ge 0.
\end{array}
\]
Let us show that $\widetilde x^t=(0,t,1)$ and $\widehat x^t=(t,0,1)$ are nondegenerate Karush-Kuhn-Tucker points of $\mathcal{S}$ if $0<t<\frac{1}{2}$. 
We consider $\widetilde x^t$ first.
Clearly, it is feasible with $\mathcal{H}(\widetilde x^t)=\{2\}$, $\mathcal{N}_1\left(\widetilde x^t\right)=\{1,3\}$, $\mathcal{N}_2\left(\widetilde x^t\right)=\emptyset$, and it holds
\[
\begin{pmatrix}
    -1+2t\\-1\\2
\end{pmatrix}
=
-\widetilde \eta^t_{2}\begin{pmatrix}
    1\\1\\-t
\end{pmatrix}
+\widetilde \nu^t_{1,1}\begin{pmatrix}
    1\\0\\0
\end{pmatrix}
+\widetilde \nu^t_{1,3}\begin{pmatrix}
    0\\0\\1
\end{pmatrix}
\]
with the unique multipliers $\widetilde \eta^t_{2}=1$, $\widetilde \nu^t_{1,1}=2t$, and $\widetilde \nu^t_{1,3}=2-t$.
Thus, we conclude that $\widetilde x^t$ is a Karush-Kuhn-Tucker of $\mathcal{S}$ fulfilling ND1 and ND2.
Moreover, %we have:
%\[
%\mathcal{T}^{\mathcal{S}}_{\tilde x^t}=\left\{\xi \in \R^3\,\left\vert \,\xi_1=\xi_2=\xi_3=0\right.\right\}.
%\]
%Thus, 
ND3 is trivially fulfilled, and for the quadratic index of $\widetilde x^t$ we have $\widetilde q=0$.
Next, we consider $\widehat x^t$.
It is feasible with $\mathcal{H}(\widehat x^t)=\{1,2\}$, $\mathcal{N}_1\left(\widehat x^t\right)=\{3\}$, $\mathcal{N}_2\left(\widehat x^t\right)=\emptyset$. Further, it holds
\[
\begin{pmatrix}
    -1\\-1+2t\\2
\end{pmatrix}
=
-\widehat \eta^t_{1}\begin{pmatrix}
    1\\-2t\\t
\end{pmatrix}
-\widehat \eta^t_{2}\begin{pmatrix}
    1\\1\\-t
\end{pmatrix}
+\widehat \nu^t_{1,3}\begin{pmatrix}
    0\\0\\1
\end{pmatrix}
\]
with the unique multipliers $\widehat \eta^t_{1}=\frac{2t}{1+2t}$, $\widehat \eta^t_{2}=\frac{1}{1+2t}$, and $\widehat \nu^t_{1,3}=2+\frac{2t^2-t}{1+2t}$.
Consequently, $\widehat x^t$ is a Karush-Kuhn-Tucker point of $\mathcal{S}$ fulfilling ND1 and ND2.
In addition, %we have:
%\[
%\mathcal{T}^{\mathcal{S}}_{\hat x^t}=\left\{\xi \in \R^3\,\left\vert \,\xi_1=\xi_2=\xi_3=0\right.\right\}.
%\]
%Hence, 
ND3 is again trivially fulfilled at $\widehat x^t$, and for its quadratic index we have $\widehat q=0$.
Note that the sequences $\widetilde x^t$ and $\widehat x^t$ converge to the same limiting point $\bar x=(0,0,1)$ for $t \to 0$. From Example \ref{ex:NDC4} we know that the latter is a nondegenerate minimizer of the underlying MPCC. 
As mentioned above, it thus satisfies MPCC-LICQ, NDC2, and MPCC-SOC.
\color{black}
Hence, Theorem \ref{thm:still} can be applied and a unique local minimizer of $\mathcal{S}^=$ in a neighborhood of $\bar x$ can be guaranteed. 
Note that it can be neither $\widetilde x^t$ nor $\widehat x^t$, since they are infeasible for $\mathcal{S}^=$.
\color{black}
In contrast to this, a neighborhood of the minimizer $\bar x$ of MPCC
contains two different nondegenerate minimizers of $\mathcal{S}$.
%, namely $\widetilde x^t$ and $\widehat x^t$. 
This is due to the violation of NDC4 at $\bar x$.
\qed
\end{example}

Example \ref{ex:2min} demonstrates that there might be multiple nondegenerate minimizers of the Scholtes regularization converging to one and the same nondegenerate minimizer of MPCC. This instability is caused by a possible violation of NDC4 and can be even worse. 
The following Example \ref{ex:continuum} demonstrates that for a nondegenerate saddle point of MPCC the number of converging sequences of (degenerate) Karush-Kuhn-Tucker points of $\mathcal{S}$ might be even infinite.

%%%%%%%%%%%%%%%%%%%%%%%%%%%%%%%%%%%%%%
%Example Continuum of KKT
%%%%%%%%%%%%%%%%%%%%%%%%%%%%%%%%%%%%%%

%%%%%%%%%%%% Anomaly Degeneracy

\begin{example}[Continuum of Karush-Kuhn-Tucker points]
\label{ex:continuum}
    We consider the following MPCC with
$n=4$ and $\kappa=3$: 
\[
\begin{array}{ll}
\min\limits_x & -x_1+x_3^2+x_4^2\\
\text{s.t.}& x_1\cdot (2-x_3) = 0, x_1\ge 0, 2-x_3 \ge 0,\\
& (x_1+x_2)\cdot (1-x_4) = 0, x_1+x_2\ge 0, 1-x_4 \ge 0,\\
& (x_3-1)\cdot (1+x_4) = 0, x_3-1\ge 0, 1+x_4 \ge 0.
\end{array}
\]
We show that $\bar x=(0,0,1,0)$ is a nondegenerate C-stationary point. Indeed, it holds
\[
\begin{pmatrix}
    -1\\0\\2\\0
\end{pmatrix}
=
\bar \sigma_{1,1}\begin{pmatrix}
    1\\0\\0\\0
\end{pmatrix}
+
\bar \sigma_{1,2}\begin{pmatrix}
    1\\1\\0\\0
\end{pmatrix}
+
\bar \sigma_{1,3}\begin{pmatrix}
    0\\0\\1\\0
\end{pmatrix}
\]
with the unique multipliers $\bar \sigma_{1,1}=-1$, $\bar \sigma_{1,2}=0$ and $\bar \sigma_{1,3}=2$.
%Thus, we conclude that $\bar x$ is a C-stationary point.
It is straightforward to verify that it is nondegenerate, i.e. NDC1-NDC3 hold, but NDC4 is violated.
% but it does not fulfill MPCC-SSOSC. 
%Since we have for its C-index $c=0$ it is also a minimizer of the considered MPCC. 
%
Let us focus on the corresponding Scholtes regularization 
\[
\begin{array}{rll}
\mathcal{S}:\min\limits_x & -x_1+x_3^2+x_4^2\\
\text{s.t.}& x_1\cdot (2-x_3) \le t, x_1\ge 0, 2-x_3 \ge 0,\\
& (x_1+x_2)\cdot (1-x_4) \le t, x_1+x_2\ge 0, 1-x_4 \ge 0,\\
& (x_3-1)\cdot (1+x_4) \le t, x_3-1\ge 0, 1+x_4 \ge 0.
\end{array}
\]
We claim that $x_z=(t,z,1,0)$ is a Karush-Kuhn-Tucker point of $\mathcal{S}$ for any $z \in [-t,0]$ with a fixed $t\leq 2$. We check for $z=-t$:
\[
\begin{pmatrix}
    -1\\0\\2\\0
\end{pmatrix}
=
-\eta_1\begin{pmatrix}
    1\\0\\-t\\0
\end{pmatrix}
+
\nu_{1,2}\begin{pmatrix}
    1\\1\\0\\0
\end{pmatrix}
+
\nu_{1,3}\begin{pmatrix}
    0\\0\\1\\0
\end{pmatrix}
\]
with the unique multipliers $\eta_1=1$, $\nu_{1,2}=0$, and $\nu_{1,3}=2-t$. Here, ND2 is violated.
Likewise, we get for $z=0$:
\[
\begin{pmatrix}
    -1\\0\\2\\0
\end{pmatrix}
=
-\eta_1\begin{pmatrix}
    1\\0\\-t\\0
\end{pmatrix}
-\eta_2\begin{pmatrix}
    1\\1\\0\\-t
\end{pmatrix}
+
\nu_{1,3}\begin{pmatrix}
    0\\0\\1\\0
\end{pmatrix}
\]
with the unique multipliers $\eta_1=1$, $\eta_2=0$, and $\nu_{1,3}=2-t$. ND2 is violated also here. Finally,
for $z\in(-t,0)$:
\[
\begin{pmatrix}
    -1\\0\\2\\0
\end{pmatrix}
=
-\eta_1\begin{pmatrix}
    1\\0\\-t\\0
\end{pmatrix}
+
\nu_{1,3}\begin{pmatrix}
    0\\0\\1\\0
\end{pmatrix}
\]
with the unique multipliers $\eta_1=1$ and $\nu_{1,3}=2-t$. Here, ND3 is violated.
Overall, we conclude that $x_z$ is degenerate for any $z\in [-t,0]$.
\color{black}
We additionally point out that all $x_z$ are infeasible for $\mathcal{S}^=$.
\color{black}
\qed
\end{example}

%%%%%%%%%%%% Anomaly Degeneracy

Finally, Example \ref{ex:continuum2} shows that infinitely many minimizers of the Scholtes regularization might occur in a neighborhood of a C-stationary point of MPCC, even if it is a nondegenerate minimizer. 
%This anomaly is once again caused by the violation of NDC4.

\begin{example}[Continuum of minimizers]
\label{ex:continuum2}
    We consider the following MPCC with
$n=2$ and $\kappa=1$: 
\[
\begin{array}{ll}
\min\limits_x & (x_1-1)^2\\
\text{s.t.}& x_1\cdot x_2 = 0, x_1\ge 0, x_2 \ge 0.
\end{array}
\]
We show that $\bar x=(1,0)$ is a nondegenerate C-stationary point. Indeed, it holds
\[
\begin{pmatrix}
    0\\0
\end{pmatrix}
=
\bar \sigma_{2,1}\begin{pmatrix}
    0\\1
\end{pmatrix}
\]
with the unique multiplier $\bar \sigma_{2,1}=0$.
%Thus, we conclude that $\bar x$ is a C-stationary point.
It is straightforward to verify that it is nondegenerate, i.e. NDC1-NDC3 hold, but NDC4 is violated. It does not surprise that any neighborhood of the minimizer $\bar x$ of MPCC
contains infinitely many minimizers of the corresponding  
% but it does not fulfill MPCC-SSOSC. 
%Since we have for its C-index $c=0$ it is also a minimizer of the considered MPCC. 
 Scholtes regularization 
\[
\begin{array}{rll}
\mathcal{S}:\min\limits_x & (x_1-1)^2\\
\text{s.t.}& x_1\cdot x_2 \le t, x_1\ge 0, x_2 \ge 0.
\end{array}
\]
Indeed, $x_z=(1,z)$ are minimizers of 
$\mathcal{S}$ for any $z \in [0,t]$. 
\color{black}
Note that $\bar x$ satisfies MPCC-LICQ, NDC2, and MPCC-SOC. Hence, Theorem \ref{thm:still} can be applied again and the unique minimizer of $\mathcal{S}^=$ is $x_t$.
\color{black}
%The latter follows easily, since the gradient of the objective function always vanishes at $x_z$. 
%\[
%\nabla f\left(x_z\right)=\begin{pmatrix}
%    0\\0
%\end{pmatrix}.
%\]
\begin{comment}
We check for $z=0$:
\[
\begin{pmatrix}
    0\\0
\end{pmatrix}
=
\nu_{2,1}\begin{pmatrix}
    0\\1
\end{pmatrix}
\]
with the unique multiplier $\nu_{2,1}=0$.
Likewise, we get for $z=t$:
\[
\begin{pmatrix}
    0\\0
\end{pmatrix}
=
-\eta_1\begin{pmatrix}
    t\\1
\end{pmatrix}
\]
with the unique multipliers $\eta_1=0$. Finally,
for $z\in(0,t)$:
\[
\begin{pmatrix}
    0\\0
\end{pmatrix}
=
\begin{pmatrix}
    0\\0
\end{pmatrix}.
\]
\end{comment}
%Note that $x_z$ is degenerate for any $z\in [0,t]$.
%We also note that all $x_z$ obviously solve $\mathcal{S}$. 
\qed
\end{example}

%%%%%%%%%%% Results for S^=

\color{black}
Due to the structural differences of $\mathcal{S}$ and $\mathcal{S}^=$, NDC4 is not relevant in the context of the latter, as observed in Examples \ref{ex:2min}-\ref{ex:continuum2}. In particular, using the same arguments as in the respective proofs, it is possible to deduce similar results for $\mathcal{S}^=$ as in Theorems \ref{thm:kktsequence} and \ref{thm:wellposedness} without relying on NDC4. Thus, nondegeneracy and the topological type locally prevail under just NDC1-NDC3.

   \begin{theorem}[Convergence for $\mathcal{S}^=$]
\label{thm:kktstill}
%Let Assumption \ref{ass:c} be fulfilled and $\varepsilon \le \frac{1}{n-s}$. 
Suppose a sequence of nondegenerate Karush-Kuhn-Tucker points $x^{t} \in M^{\mathcal{S}^=}$ of 
$\mathcal{S}^=$ with quadratic index $q$ converges to $\bar x$ for $t \to 0$. If $\bar x \in M$ is a nondegenerate C-stationary point of MPCC with multipliers $(\bar \sigma, \bar \varrho)$, then we have for its C-index $CI= q$.
\end{theorem}

We point out that the results of next Theorem \ref{thm:wpstill}, except of tracing the corresponding indices, have been stated in Remark 5.1 from \cite{still:2007}.

\begin{theorem}[Well-posedness for $\mathcal{S}^=$]
\label{thm:wpstill}
   Let $\bar x \in M$ be a nondegenerate C-stationary point of MPCC with C-index $c$. Then, for all sufficiently small $t$ there exists a nondegenerate Karush-Kuhn-Tucker point $x^t \in M^{\mathcal{S}^=}$ of $\mathcal{S}^=$ within a neighborhood of $\bar x$, which has the same quadratic index $c$.
   Moreover, for any fixed $t$ sufficiently small, such $x^t$ is the unique Karush-Kuhn-Tucker point of $\mathcal{S}^=$ in a sufficiently small neighborhood of $\bar x$.
\end{theorem} 

\color{black}

\section*{Conclusion}
 We successfully related C-stationary points of MPCC with Karush-Kuhn-Tucker points of its Scholtes regularization in regard to their nondegeneracy, uniqueness, and index. 
This enables us to study the global structure of the Scholtes regularization for MPCC next. In particular, the question on the homotopy equivalence of the corresponding lower level sets of MPCC and its Scholtes regularization can be addressed now. Another interesting line of future research is to apply the described approach to other regularizations for MPCC known from the literature. By doing so, we expect to identify those regularizations which guarantee the preservation of the above invariants under less restrictive additional conditions for MPCC. Needless to say that also regularizations for mathematical programs with other nonsmooth constraints, such as vanishing, switching, orthogonality-type, disjunctive etc., can be analogously studied. Overall, we conclude that this paper opens a door for a better understanding of the structural aspects of regularizations in nonsmooth optimization.

%%%%%%%%%%%%%%%%%%%%%%%%%%%%%%%%%%%%%%%
%Appendix
%%%%%%%%%%%%%%%%%%%%%%%%%%%%%%%%%%%%%%%

%\newpage
\color{black}
\section*{Appendix}
%\begin{proof}
\color{black}
\subsection*{Proof of Theorem \ref{thm:kktsequence}.}
Throughout the proof, we assume for the sake of simplicity and without loss of generality
that $a_{10}\left(\bar x\right)=\emptyset$.
Further, we set $p=\left\vert a_{00}\left(\bar x\right)\right\vert$ and assume that the considered problems are given in standard form, analogously to \cite{still:2007}:
\begin{equation}
\label{eq:mpccstd}
\begin{array}{lrll}
   \text{MPCC}: &\min\limits_x& f(x)&\\
    &\text{s.t.}
    &x_j\cdot x_{\kappa + j}=0&j\in\{1,\ldots,p\},\\
    &&x_j\cdot F_{2,j}(x)=0&j\in\{p+1,\ldots,\kappa\},\\
    &&x_j\ge 0& j\in\{1,\ldots,\kappa+p\},\\
    &&F_{2,j}(x)\ge 0&j\in \{p+1,\ldots,\kappa\},
\end{array}
\end{equation}

and
\begin{equation}
\label{eq:scholtesstd}
    \begin{array}{lrll}
   \mathcal{S}: &\min\limits_x& f(x)&\\
    &\text{s.t.}
    &x_j\cdot x_{\kappa + j} \le t&j\in\{1,\ldots,p\},\\
    &&x_j\cdot F_{2,j}(x)\le t&j\in\{p+1,\ldots,\kappa\},\\
    &&x_j\ge 0& j\in\{1,\ldots,\kappa+p\},\\
    &&F_{2,j}(x)\ge 0&j\in \{p+1,\ldots,\kappa\}.
\end{array}
\end{equation}

\begin{comment}
\[
f(x)=\sum\limits_{j=p+1}^{\kappa}\bar \sigma_{1,j}x_j+\sum\limits_{j=1}^p\left(\bar \varrho_{1,j}x_j+\bar \varrho_{2,j}x_{\kappa+j}\right)+q(x) \mbox{ with } Dq(\bar x)=0.
\]
\end{comment}
We recall that for all $x$ sufficiently close to $\bar x$ it holds $F_{2,j}(x)>0$, $j\in\{p+1,\ldots,\kappa\}$. Hence, the Hessians of the respective Lagrange functions are
\[
D^2L(x)=D^2f(x),
\]
%D^2F-Teil ergänzen
\[
\begin{array}{rcl}

D^2L^{\mathcal{S}}(x)&=&D^2f(x)+\displaystyle\sum\limits_{j=1}^p\eta_j\left(e_j e_{\kappa+j}^T+e_{\kappa+j} e_{j}^T\right)
+\displaystyle\sum\limits_{j=p+1}^\kappa \eta_j\left(e_j D F_{2,j}\left(x\right)+D^T F_{2,j}\left(x\right) e_j^T\right).
\end{array}
\]
Further, the tangent space at the C-stationary point $\bar x$ of MPCC calculates as
\[
\mathcal{T}_{\bar x}=\left\{\xi\in \R^n\,\left\vert\,\xi_j=0, j\in\{1,\ldots,\kappa+p\}\right.\right\}.
\]
For the tangent space at the Karush-Kuhn-Tucker point $x$ of $\mathcal{S}$ we get
\[
    \mathcal{T}^\mathcal{S}_{x}=\left\{
\xi \in \R^{n}\,\left\vert\, \begin{array}{l}
\left( e_j^T \cdot x_{\kappa+j}+x_j\cdot e^T_{\kappa+j}\right)\xi=0, j \in \mathcal{H}\left(x\right)\cap\{1,\ldots,p\},\\
\left( e_j^T \cdot F_{2,j}\left(x\right)+x_j\cdot D F_{2,j}\left(x\right)\right)\xi=0, j \in \mathcal{H}\left(x\right)\cap\{p+1,\ldots,\kappa\},\\
\xi_j=0, j \in \mathcal{N}_1\left(x\right),
\xi_{\kappa+j}=0, j \in \mathcal{N}_2\left(x\right)
\end{array}
\right.\right\}.
\]
The remaining proof will be divided into five major steps.

{\bf Step 1.}
We denote the codimensions of $\mathcal{T}_{\bar x}$ and $\mathcal{T}^\mathcal{S}_{x^t}$, respectively, as 
\[\alpha=\left\vert a_{01}(\bar x)\right\vert+2\left\vert a_{00}(\bar x)\right\vert, \quad 
\alpha^\mathcal{S}_t=
\left\vert \mathcal{H}\left(x^t\right)\right\vert +
\left\vert  \mathcal{N}_1\left(x^t\right)\right\vert +
\left\vert  \mathcal{N}_2\left(x^t\right)\right\vert.
\]
We claim that 
\begin{equation}
\label{eq:step1}
\begin{array}{rcl}
\alpha - \alpha_t^\mathcal{S} &\leq&
\left\vert  a_{01}^0\left(\bar x\right)\right\vert +
%\left\vert  a_{10}^0\left(\bar x\right)\right\vert +
\left\vert  a_{00}^-\left(\bar x\right)\right\vert. 
\end{array}
\end{equation}
%
\begin{comment}
We consider
\[
\begin{array}{rcl}
\alpha-\alpha^\mathcal{S}&=&
\left\vert  a_{01}\left(\bar x\right)\right\vert +
\left\vert  a_{10}\left(\bar x\right)\right\vert +
2\left\vert  a_{00}\left(\bar x\right)\right\vert\\&&-
\left\vert \mathcal{H}\left(x^t\right)\right\vert-
\left\vert  \mathcal{N}_1\left(x^t\right)\right\vert -
\left\vert  \mathcal{N}_2\left(x^t\right)\right\vert
\end{array}
\]
\end{comment}
To see this, we first recall that 
$a_{01}\left(\bar x\right)%, a_{10}\left(\bar x\right)
$, 
and $a_{00}\left(\bar x\right)$ are disjoint.
Hence, it holds
\[
\begin{array}{rcl}
\alpha^\mathcal{S}_t&=&
\left\vert \mathcal{H}\left(x^t\right)\cap a_{01}\left(\bar x\right)\right\vert + 
%\left\vert \mathcal{H}\left(x^t\right)\cap a_{10}\left(\bar x\right)\right\vert +
\left\vert \mathcal{H}\left(x^t\right)\cap a_{00}\left(\bar x\right)\right\vert+
\left\vert  \mathcal{N}_1\left(x^t\right)\cap a_{01}\left(\bar x\right)\right\vert +
%\left\vert  \mathcal{N}_1\left(x^t\right)\cap a_{10}\left(\bar x\right)\right\vert +
\left\vert  \mathcal{N}_1\left(x^t\right)\cap a_{00}\left(\bar x\right)\right\vert \\&&+
\left\vert  \mathcal{N}_2\left(x^t\right)\cap a_{01}\left(\bar x\right)\right\vert +
%\left\vert  \mathcal{N}_2\left(x^t\right)\cap a_{10}\left(\bar x\right)\right\vert +
\left\vert  \mathcal{N}_2\left(x^t\right)\cap a_{00}\left(\bar x\right)\right\vert.
\end{array}
\]
By continuity arguments, we have for $t$ sufficiently small:
\[
%\mathcal{N}_1\left(x^t\right)\cap a_{10}\left(\bar x\right)=
\mathcal{N}_2\left(x^t\right)\cap a_{01}\left(\bar x\right)=\emptyset.
\]
From $\mathcal{H}\left(x^t\right) \cap \mathcal{N}_1\left(x^t\right)=\mathcal{H}\left(x^t\right) \cap \mathcal{N}_2\left(x^t\right) = \emptyset$ -- together with Corollary \ref{cor:activeindex:scholtes}a)--c) -- we have:
\[
   \mathcal{H}\left(x^t\right)\cap a_{01}^+\left(\bar x\right) = %\mathcal{H}\left(x^t\right)\cap a_{10}^+\left(\bar x\right) = 
   \mathcal{H}\left(x^t\right)\cap a_{00}^+\left(\bar x\right)= \emptyset.  
\]
and 
\[
   \mathcal{N}_1\left(x^t\right)\cap a_{01}^-\left(\bar x\right) = %\mathcal{N}_2\left(x^t\right)\cap a_{10}^-\left(\bar x\right) =
   \mathcal{N}_1\left(x^t\right)\cap a_{00}^-\left(\bar x\right)= \mathcal{N}_2\left(x^t\right)\cap a_{00}^-\left(\bar x\right)=\emptyset.  
\]
Recalling $a_{00}^0\left(\bar x\right)= \emptyset$, due to NDC2, we then get:  
\[
\begin{array}{rcl}
\alpha^\mathcal{S}_t&=&
\left\vert \mathcal{H}\left(x^t\right)\cap a_{01}^-\left(\bar x\right)\right\vert +
\left\vert \mathcal{H}\left(x^t\right)\cap a_{01}^0\left(\bar x\right)\right\vert
+
%\left\vert \mathcal{H}\left(x^t\right)\cap a_{10}^-\left(\bar x\right)\right\vert %+\left\vert \mathcal{H}\left(x^t\right)\cap a_{10}^0\left(\bar x\right)\right\vert\\&&+
\left\vert \mathcal{H}\left(x^t\right)\cap a_{00}^-\left(\bar x\right)\right\vert \\&&+
\left\vert  \mathcal{N}_1\left(x^t\right)\cap a_{01}^0\left(\bar x\right)\right\vert +
\left\vert  \mathcal{N}_1\left(x^t\right)\cap a_{01}^+\left(\bar x\right)\right\vert
+\left\vert  \mathcal{N}_1\left(x^t\right)\cap a_{00}^+\left(\bar x\right)\right\vert 
%+\left\vert  \mathcal{N}_2\left(x^t\right)\cap a_{10}^0\left(\bar x\right)\right\vert +
%\left\vert  \mathcal{N}_2\left(x^t\right)\cap a_{10}^+\left(\bar x\right)\right\vert
+\left\vert  \mathcal{N}_2\left(x^t\right)\cap a_{00}^+\left(\bar x\right)\right\vert\\
&=&
\left\vert a_{01}^-\left(\bar x\right)\right\vert +
\left\vert \mathcal{H}\left(x^t\right)\cap a_{01}^0\left(\bar x\right)\right\vert
+
%\left\vert a_{10}^-\left(\bar x\right)\right\vert +
%\left\vert \mathcal{H}\left(x^t\right)\cap a_{10}^0\left(\bar x\right)\right\vert +
\left\vert  a_{00}^-\left(\bar x\right)\right\vert \\&&+
\left\vert  \mathcal{N}_1\left(x^t\right)\cap a_{01}^0\left(\bar x\right)\right\vert +
\left\vert  a_{01}^+\left(\bar x\right)\right\vert
+\left\vert  a_{00}^+\left(\bar x\right)\right\vert %\\&& 
%+\left\vert  \mathcal{N}_2\left(x^t\right)\cap a_{10}^0\left(\bar x\right)\right\vert +
%\left\vert  a_{10}^+\left(\bar x\right)\right\vert
+\left\vert a_{00}^+\left(\bar x\right)\right\vert
\\ 
&=& \left\vert  a_{01}^-\left(\bar x\right)\right\vert +
%\left\vert  a_{10}^-\left(\bar x\right)\right\vert +
\left\vert  a_{00}^-\left(\bar x\right)\right\vert +
\left\vert  a_{01}^+\left(\bar x\right)\right\vert + 
%\left\vert  a_{10}^+\left(\bar x\right)\right\vert +
2\left\vert a_{00}^+\left(\bar x\right)\right\vert 
 %\\&& 
 +\left\vert \left(\mathcal{H}\left(x^t\right)\cup \mathcal{N}_1\left(x^t\right)\right) \cap a_{01}^0\left(\bar x\right)\right\vert.
 %+ \left\vert \left(\mathcal{H}\left(x^t\right)\cup \mathcal{N}_2\left(x^t\right)\right) \cap a_{10}^0\left(\bar x\right)\right\vert
\end{array}
\]
\textcolor{black}{Altogether, we obtain with the above:}
\[
\begin{array}{rcl}
\alpha - \alpha_t^\mathcal{S} &=&
\left\vert  a_{01}^0\left(\bar x\right)\right\vert +
%\left\vert  a_{10}^0\left(\bar x\right)\right\vert +
\left\vert  a_{00}^-\left(\bar x\right)\right\vert 
%\\&& 
-\left\vert \left(\mathcal{H}\left(x^t\right)\cup \mathcal{N}_1\left(x^t\right)\right) \cap a_{01}^0\left(\bar x\right)\right\vert.
%- \left\vert \left(\mathcal{H}\left(x^t\right)\cup \mathcal{N}_2\left(x^t\right)\right) \cap a_{10}^0\left(\bar x\right)\right\vert.
\begin{comment}
\left\vert \left\{j\in a_{01}\left(\bar x\right)\, \left\vert \, \sigma_{1,j}=0, j \not\in \mathcal{H}\left(x^t\right)\cup \mathcal{N}_1\left(x^t\right)\right.\right\}\right\vert\\
&&+\left\vert \left\{j\in a_{10}\left(\bar x\right)\, \left\vert \, \sigma_{2,j}=0, j \not\in \mathcal{H}\left(x^t\right)\cup \mathcal{N}_2\left(x^t\right)\right.\right\}\right\vert\\
&&+\left\vert  a_{00}^-\left(\bar x\right)\right\vert
\\
%&\le&
%\left\vert\left\{j\in a_{01}\left(\bar x\right)\, \left\vert \, \sigma_{1,j}=0\right.\right\}\right\vert
%+\left\vert\left\{\in a_{10}\left(\bar x\right)\, \left\vert \, \sigma_{2,j}=0\right.\right\}\right\vert
%+\left\vert  a_{00}^-\left(\bar x\right)\right\vert\\.
\end{comment}
\end{array}
\]
From here, the claimed inequality (\ref{eq:step1}) follows.

{\bf Step 2.}
Let $\mathcal{T}\subset \R^{n}$ be a linear subspace. We denote the number of negative
eigenvalues of $D^2 L^\mathcal{S}\left( x^t\right)\restriction_{\mathcal{T}}$ by $QI^\mathcal{S}_{t,\mathcal{T}}$. Analogously, $QI_{t,\mathcal{T}}$ stands for the number of negative eigenvalues of $D^2 L\left( x^t\right)\restriction_{\mathcal{T}}$ and $\overline{QI}_{\mathcal{T}}$ stands for the number of negative eigenvalues of $D^2 L\left(\bar x \right)\restriction_{\mathcal{T}}$. 
We have the following relation between the involved Hessians of the Lagrange functions

\begin{equation}
\label{eq:Hessetandbar}
\begin{array}{rcl}
\displaystyle
D^2 L^\mathcal{S}\left( x^t\right) &=& 
\displaystyle D^2 L\left( x^t\right) + \sum\limits_{j=1}^p\eta_j^t\left(e_j\cdot e_{\kappa+j}^T+e_{\kappa+j}\cdot e_{j}^T\right)
\\&&+\displaystyle\sum\limits_{j=p+1}^\kappa \eta_j^t\left(e_j\cdot D F_{2,j}\left(x^t\right)+D^T F_{2,j}\left(x^t\right)\cdot e_j^T\right).
\end{array}
\end{equation}
\begin{comment}
Suppose $\xi \in \mathcal{T}_{\bar x}$. We consider
\[
\xi^T \left( \displaystyle\sum\limits_{j=1}^{\kappa} 
\eta_j^{t} \left(F_{2,j}\left(x^t\right)D^2F_{1,j}\left(x^t\right)+\nabla F_{2,j}\left(x^t\right)DF_{1,j}\left(x^t\right)
+\nabla F_{1,j}\left(x^t\right)DF_{2,j}\left(x^t\right)+F_{1,j}\left(x^t\right)D^2F_{2,j}\left(x^t\right)\right)\right)\xi.
\]
 Note, that $\eta_i=0$ for all $i \notin \mathcal{H}\left(\bar x^t\right)$, due to (\ref{eq:kkt-2}).
 %$i\in a_{01}^+\left(bar x\right) \cup a_{10}^+\left(bar x\right) \cup a_{00}^+\left(bar x\right)$, due to Corollary \ref{cor:activeindex:scholtes}.
\end{comment}
We claim that for all $t$ sufficiently small it holds: 
\begin{equation}
\label{eq:step2a}
QI^\mathcal{S}_{t,\mathcal{T}_{\bar x}}=\overline{QI}_{\mathcal{T}_{\bar x}}.
\end{equation}
Clearly, we have for $\xi\in \mathcal{T}_{\bar x}$ that 
%e_j^T\xi =0$ for $j \in a_{01}\left(\bar x\right)=\{p+1,\ldots,\kappa\}$. Moreover, for $j \in a_{00}\left(\bar x\right)=\{1,\ldots,p\}$ we also have $e_j^T\xi =0$ and $e_{\kappa+j}^T\xi =0$.
\[
\xi^T \left(\sum\limits_{j=1}^p\eta^t_j\left(e_j\cdot e_{\kappa+j}^T+e_{\kappa+j}\cdot e_{j}^T\right)+\sum\limits_{j=p+1}^\kappa \eta^t_j\left(e_j\cdot D F_{2,j}\left(x^t\right)+D^T F_{2,j}\left(x^t\right)\cdot e_j^T\right)
\right)\xi = 0.
\]
Hence, it holds:
\begin{equation}
    \label{eq:step2limits}
\xi^T D^2 L^\mathcal{S}\left( x^t\right) \xi \rightarrow \xi^T D^2 L\left( \bar x\right) \xi \mbox{ for } t \rightarrow 0.
\end{equation}
In view of NDC3, the formula (\ref{eq:step2a}) is therefore valid.

{\bf Step 3a.}
We claim that the numbers of positive and negative eigenvalues of $D^2 L^\mathcal{S}\left( x^t\right)\restriction_{ \mathcal{T}_{\bar x}}$ and of $D^2 L^\mathcal{S}\left( x^t\right)\restriction_{\mathcal{T}^t}$, respectively, coincide, where
\[
  \mathcal{T}^t =\left\{
\xi \in \R^{n}\,\left\vert\, \begin{array}{l}
\left( e_j^T \cdot F_{2,j}\left(x^t\right)+x_j\cdot D F_{2,j}\left(x^t\right)\right)\xi=0,
j \in a_{01}^-\left(\bar x\right)%\cup a_{10}^-\left(\bar x\right)
,\\
\xi_j=0, j \in a_{01}^+\left(\bar x\right) \cup a_{01}^0\left(\bar x\right) \cup a_{00}\left(\bar x\right),
\xi_{\kappa+j}=0, j \in %a_{10}^+\left(\bar x\right) \cup a_{10}^0\left(\bar x\right) \cup 
a_{00}\left(\bar x\right)
\end{array}
\right.\right\}.
\]
To see this, let $\left\{\lambda^+_1,\ldots,\lambda^+_{k^+}\right\}$ be the positive 
eigenvalues of 
$D^2 L^\mathcal{S}\left( x^t\right)\restriction_{ \mathcal{T}_{\bar x}}$ with corresponding eigenvectors 
$\left\{\xi^+_1,\ldots,\xi^+_{k^+}\right\} \subset \mathcal{T}_{\bar x}$.
%Hence, for all $k=1,\ldots, k^+$:
%\[
%{\xi^+_k}^T D^2 L^\mathcal{S}\left( x^t\right)\xi^+_k>0.
%\]
We further equivalently rewrite the tangent space
\[
\mathcal{T}_{\bar x}=\left\{\xi \in \R^n \,\left\vert \,
\begin{array}{l}
\left( e_j^T \cdot F_{2,j}\left(\bar x\right)+\bar x_j\cdot D F_{2,j}\left(\bar x\right)\right)\xi=0,
j \in a_{01}^-\left(\bar x\right)%\cup a_{10}^-\left(\bar x\right)
,\\
\xi_j=0, j \in a_{01}^+\left(\bar x\right) \cup a_{01}^0\left(\bar x\right) \cup a_{00}\left(\bar x\right),
\xi_{k+j}=0, j \in  %a_{10}^+\left(\bar x\right) \cup a_{10}^0\left(\bar x\right) \cup 
a_{00}\left(\bar x \right)
\end{array}
\right.\right\}.
\]
In view of MPCC-LICQ, the implicit function theorem can be applied. It provides the existence of
$\delta_1,\delta_2>0$, such that for all $k=1,\ldots,k^+$ and $t<\delta_1$ there exists $\xi_{k,t}$ with $\left\| \xi_{k,t}-\xi^+_k\right\|<\delta_2$ and $\xi_{k,t}\in \mathcal{T}^t$.
 We can choose $t$ even smaller, so that $\xi_{1,t},\ldots,\xi_{k^+,t}$ remain linearly independent and for all $k=1,\ldots,k^+$ it holds:
\[
{ \xi_{k,t}}^T D^2 L^\mathcal{S}\left( x^t\right) \xi_{k,t}>0.
\]
Hence, $D^2 L^\mathcal{S}\left( x^t\right)\restriction_{\mathcal{T}^t}$ has at least $k^+$ positive eigenvalues. If we repeat the above reasoning for negative eigenvalues, the matrix $D^2 L^\mathcal{S}\left( x^t\right)\restriction_{\mathcal{T}^t}$ has at least as many negative eigenvalues as $D^2 L^\mathcal{S}\left( x^t\right)\restriction_{\mathcal{T}_{\bar x}}$. Additionally, we see that the dimensions of $\mathcal{T}_{\bar x}$ and $\mathcal{T}^t$ coincide. Moreover, the matrix $D^2 L^\mathcal{S}\left( x^t\right)\restriction_{\mathcal{T}_{\bar x}}$ is regular in view of NDC3 and Step 2. Altogether, the assertion follows.

{\bf Step 3b.}
We claim that for any $t$ sufficiently small we have 
\begin{equation}
    \label{eq:step3b}
    \mathcal{T}^t\subset \mathcal{T}^\mathcal{S}_{x^t}.
\end{equation} 
\begin{comment}
We claim 
\[
\mbox{dim}\left(\mathcal{T}^\mathcal{S}_{x^t}\right)=\mbox{dim}\left(\mathcal{T}^{\prime}\right)-\left|a_{01}^0\left(\bar x \right)\right\backslash\left(\mathcal{H}\left(x^t\right)\cup \mathcal{N}_1\left(x^t\right)\right)|-\left|a_{10}^0\left(\bar x \right)\backslash\left(\mathcal{H}\left(x^t\right)\cup \mathcal{N}_2\left(x^t\right)\right)\right|-\left|a_{00}^-\left(\bar x \right)\right|
\]
\end{comment}
In fact, $\mathcal{H}\left(x^t\right)$ and $\mathcal{N}_1\left(x^t\right) \cup \mathcal{N}_2\left(x^t\right)$ are obviously disjoint. By continuity arguments, we get $a_{01}\left(\bar x\right)\cap \mathcal{N}_2\left(x^t\right)=\emptyset$ 
%and 
%$a_{10}\left(\bar x\right)\cap \mathcal{N}_1\left(x^t\right)=\emptyset$ 
for all $t$ sufficiently small. 
By combining this with Corollary \ref{cor:activeindex:scholtes}, NDC2, and $a_{10}\left(\bar x\right)=\emptyset$, we have for all $t$ sufficiently small:

\[
 \mathcal{H}\left(x^t\right)
    =a_{01}^-\left(\bar x\right)%\cup  a_{10}^-\left(\bar x\right)
    \cup  a_{00}^-\left(\bar x\right)\cup 
    \left(a_{01}^0\left(\bar x\right)\cap \mathcal{H}\left(x^t\right)\right),
    %\cup \left(a_{10}^0\left(\bar x\right)\cap \mathcal{H}\left(x^t\right)\right), 
\]
\[
       \mathcal{N}_1\left(x^t\right)  =a_{01}^+\left(\bar x\right) \cup   a_{00}^+\left(\bar x\right)
    \cup \left(a_{01}^0\left(\bar x\right)\cap \mathcal{N}_1\left(x^t\right)\right), \quad
            %a_{10}^+\left(\bar x\right) \cup 
    %\left(a_{10}^0\left(\bar x\right)\cap \mathcal{N}_2\left(x^t\right)\right) \cup
    \mathcal{N}_2\left(x^t\right)= a_{00}^+\left(\bar x\right).
\]
Thus, we get
\[
\mathcal{T}^\mathcal{S}_{x^t}=
\left\{\xi \in \R^n \,\left\vert \,
\begin{array}{l}
\left( e_j^T \cdot F_{2,j}\left(x^t\right)+x_j\cdot D F_{2,j}\left(x^t\right)\right)\xi=0,\\
j \in a_{01}^-\left(\bar x\right)%\cup  a_{10}^-\left(\bar x\right)
\cup  a_{00}^-\left(\bar x\right)\cup 
    \left(a_{01}^0\left(\bar x\right)\cap \mathcal{H}\left(x^t\right)\right)
    %\cup \left(a_{10}^0\left(\bar x\right)\cap \mathcal{H}\left(x^t\right)\right)
    ,\\
\xi_j=0, j \in a_{01}^+\left(\bar x\right) \cup   a_{00}^+\left(\bar x\right) \cup\left(a_{01}^0\left(\bar x\right)\cap \mathcal{N}_1\left(x^t\right)\right),
\xi_{\kappa+j}=0, j \in %a_{10}^+\left(\bar x\right) \cup 
    %\left(a_{10}^0\left(\bar x\right)\cap \mathcal{N}_2\left(x^t\right)\right) \cup 
    a_{00}^+\left(\bar x\right)
\end{array}
\right.\right\}.
\]
From here, (\ref{eq:step3b}) directly follows.

{\bf Step 4.} 
In view of Steps 3a) and 3b), we conclude that for all $t$ sufficiently small:
\[
QI^\mathcal{S}_{t,\mathcal{T}^\mathcal{S}_{x^t}} \leq 
\mbox{dim}  \mathcal{T}^\mathcal{S}_{x^t}
- 
\mbox{dim}  \mathcal{T}^t  +
QI^\mathcal{S}_{t,\mathcal{T}^t} 
= n-\alpha^\mathcal{S}_t-k^+
\]
as well as 
\[
QI^\mathcal{S}_{t,\mathcal{T}_{\bar x}}= \mbox{dim}  \mathcal{T}_{\bar x} - k^+=
n-\alpha-k^+.
\]
%\overset{Step\,\,2a}{=}QI^\mathcal{R}_{\left(t,\mathcal{T}^\mathcal{R}_{\left(\bar x,\bar y\right)}\right)}%=\overline{QI}^\mathcal{R}_{\left(\mathcal{T}^\mathcal{R}_{\left(\bar x,\bar y\right)}\right)}
It immediately follows that
\[
QI^\mathcal{S}_{t,\mathcal{T}^\mathcal{S}_{x^t}}\le QI^\mathcal{S}_{t,\mathcal{T}_{\bar x}}+\alpha-\alpha^\mathcal{S}_t.
\]
We further estimate:
\[
\begin{array}{rcl}
q=QI^\mathcal{S}_{t,\mathcal{T}^\mathcal{S}_{x^t}}
&\overset{(\ref{eq:step2a})}{\leq} &
\overline{QI}_{\mathcal{T}_{\bar x}}
+\alpha-\alpha^\mathcal{S}_t
\overset{(\ref{eq:step1})}{\le} \overline{QI}_{\mathcal{T}_{x}}+\left\vert  a_{01}^0\left(\bar x\right)\right\vert +
%\left\vert  a_{10}^0\left(\bar x\right)\right\vert +
\left\vert  a_{00}^-\left(\bar x\right)\right\vert\\
&=&\overline{QI}_{\mathcal{T}_{x}}+\left\vert  a_{01}^0\left(\bar x\right)\right\vert +
%\left\vert  a_{10}^0\left(\bar x\right)\right\vert + 
BI=CI+\left\vert\left\{j\in a_{01}\left(\bar x\right)\,\left\vert\,\bar \sigma_{1,j}=0\right.\right\}\right\vert %+ \left\vert\left\{j\in a_{10}\left(\bar x\right)\,\left\vert\,\bar \sigma_{2,j}=0\right.\right\} \right\vert
.
\end{array}
    \]
We thus obtain:
\[
CI \geq \max\left\{q - \left\vert\left\{j\in a_{01}\left(\bar x\right)\,\left\vert\,\bar \sigma_{1,j}=0\right.\right\}\right\vert %- \left\vert\left\{j\in a_{10}\left(\bar x\right)\,\left\vert\,\bar \sigma_{2,j}=0\right.\right\} \right\vert
, 0\right\}.
\]

{\bf Step 5a.}
We recall that for $t$ sufficiently small we have, in view of Corollary \ref{cor:activeindex:scholtes}c), $a_{00}^-\left(\bar x\right)\subset \mathcal{H}\left(x^t\right)$ and, thus, 
\[
x^t_j \not = 0, x^t_{\kappa+j} \ne 0 \mbox{ for all } j \in a_{00}^-\left(\bar x\right). 
\]
Therefore, both
$\frac{x_{\kappa+j}^t}{x_{j}^t}$ and
$\frac{x_{j}^t}{x_{\kappa+j}^t}$ are well-defined.
After considering subsequences if needed, we can assume -- without loss of generality --
that for any $j \in a_{00}^-\left(\bar x\right)$ at least one of the mentioned sequences is convergent. It is convenient to define auxiliary index sets:
\[
a_{00}^{-,1}\left(\bar x\right)
=\left\{j \in a_{00}^-\left(\bar x\right)\, \left\vert \, \frac{x^t_{\kappa+j}}{x^t_{j}} \mbox{ converges for } t \to 0 \right.\right\}
\quad \mbox{and}\quad
a_{00}^{-,2}\left(\bar x\right)
=a_{00}^{-}\left(\bar x\right) \backslash a_{00}^{-,1}\left(\bar x\right).
\]
Suppose $\bar j \in a_{00}^{-,1}\left(\bar x\right)$% and let $R_{01} \subset a_{01}^0\left(\bar x\right)$ and $R_{10} \subset a_{10}^0\left(\bar x\right)$
. We define the following two linear spaces:
\[
\mathcal{T}%_{R_{01},R_{10}}
^{\bar j}=\left\{\xi \in \R^n \,\left\vert \,
\begin{array}{l}
%\left(F_{2,j}\left(x^t\right)DF_{1,j}\left(x^t\right)+F_{1,j}\left(x^t\right)DF_{2,j}\left(x^t\right)\right)\xi=0,\\
%j \in a_{01}^-\left(\bar x\right)\cup a_{10}^-\left(\bar x\right)\cup R_{01} \cup R_{10} ,\\
\xi_j=0, j \in a_{01}\left(\bar x\right)\cup a_{00}\left(\bar x\right)\backslash\{\bar j\}%\cup \left(a_{01}^0\left(\bar x\right)\backslash R_{01}\right)
,
\xi_{\kappa+j}=0, j \in  %a_{10}\left(\bar x\right)\cup 
a_{00}\left(\bar x \right)\backslash\{\bar j\}%%\cup \left(a_{10}^0\left(\bar x\right)\backslash R_{10}\right)
,\\
 \left(\displaystyle \lim\limits_{t \rightarrow 0}\frac{x^t_{\kappa+\bar j}}{x^t_{\bar j}}e^T_{\bar j}+e_{\kappa+\bar j}^T\right)\xi=0
\end{array}
\right.\right\}
\]
and 
\[
\mathcal{T}%_{R_{01},R_{10}}
^{t,\bar j}=\left\{\xi \in \R^n \,\left\vert \,
\begin{array}{l}
%\left(F_{2,j}\left(x^t\right)DF_{1,j}\left(x^t\right)+F_{1,j}\left(x^t\right)DF_{2,j}\left(x^t\right)\right)\xi=0,\\
%j \in a_{01}^-\left(\bar x\right)\cup a_{10}^-\left(\bar x\right)\cup R_{01} \cup R_{10} ,\\
\xi_j=0, j \in a_{01}\left(\bar x\right)\cup a_{00}\left(\bar x\right)\backslash\{\bar j\}%\cup \left(a_{01}^0\left(\bar x\right)\backslash R_{01}\right)
,
\xi_{\kappa+j}=0, j \in  %a_{10}\left(\bar x\right)\cup 
a_{00}\left(\bar x \right)\backslash\{\bar j\}%\cup \left(a_{10}^0\left(\bar x\right)\backslash R_{10}\right)
,\\
 \left(\displaystyle \frac{x^t_{\kappa+\bar j}}{x^t_{\bar j}}e^T_{\bar j}+e_{\kappa+\bar j}^T\right)\xi=0
\end{array}
\right.\right\}.
\]

First, let us focus on $\mathcal{T}%_{R_{01},R_{10}}
^{\bar j}$. Clearly, we have
$\mathcal{T}_{\bar x} \subset \mathcal{T}%_{R_{01},R_{10}}
^{\bar j}$. Further, all the vectors involved in the description of $\mathcal{T}%_{R_{01},R_{10}}
^{\bar j}$ remain linearly independent as it is easy to verify.
Therefore, we have 
\[
\dim\left(\mathcal{T}^{\bar j}\right)=\dim\left(\mathcal{T}_{\bar x}\right)+1.
\]
Especially, we find $\xi^{\bar j} \in \mathcal{T}^{\bar j}$ such that
$
\xi^{\bar j}_{\bar j}\ne 0$.
%\quad DF_{2,j}\left(\bar x\right)\xi^{\bar j} = 0, j \in a_{00}\left(\bar x \right)\backslash\{\bar j\}.
Due to the linear independence of the involved vectors in $\mathcal{T}^{\bar j}$, we can apply the implicit function theorem once again. By doing so, we can find a sequence of $\xi^{t,\bar j} \in \mathcal{T}^{t,\bar j}$ with
$\xi^{t,\bar j} \to \xi^{\bar j}$ for $t \to 0$ and 
\begin{equation}
    \label{eq:Step5DF2ne0}
\xi^{t,\bar j}_{\bar j}\ne 0.
%, \quad DF_{1,j}\left(x^t\right)\xi^{t,\bar j}=DF_{2,j}\left(x^t\right)\xi^{t,\bar j}= 0, j \in a_{00}\left(\bar x \right)\backslash\{\bar j\}.
\end{equation}
Analogously, for $\bar j \in a_{00}^{-,2}\left(\bar x\right)$, 
we again find a sequence of $\xi^{t,\bar j} \in \mathcal{T}%_{R_{01},R_{10}}
^{t,\bar j}$ with
$\xi^{t,\bar j} \to \xi^{\bar j}$ for $t \to 0$ and 
\begin{equation}
    \label{eq:Step5DF1ne0}
\xi^{t,\bar j}_{\kappa+\bar j}\ne 0.
%, \quad DF_{1,j}\left(x^t\right)\xi^{t,\bar j}= DF_{2,j}\left(x^t\right)\xi^{t,\bar j}= 0, j \in a_{00}\left(\bar x \right)\backslash\{\bar j\}.
\end{equation}
We can use the same arguments as in Step 3b to conclude that $\xi^{t,\bar j} \in \mathcal{T}^{\mathcal{S}}_{x^t}$.

{\bf Step 5b.}
We claim that for any basis $\left\{\xi^1,\ldots,\xi^{\ell}\right\}$ of $\mathcal{T}^{t}$, cf.~Step 3a, the following set of vectors is linearly independent:
\[
\left\{\xi^1,\ldots,\xi^{\ell}\right\} \cup  \left\{\xi^{t,j} \,  \left\vert \, j \in a_{00}^-\left(\bar x \right) \right.\right\}.
\]
Indeed, consider a linear combination of the latter as given by
\[
\sum\limits_{i=1}^{\ell} b_i \xi^{i} + \sum\limits_{j \in a_{00}^-\left(\bar x \right)} \beta_j \xi^{t,j} = 0.
\]
For $\bar j \in a_{00}^{-,1}\left(\bar x\right)$ let us multiply both sides with $e_{\bar j}^T$ from the left.
In view of $\xi^i \in \mathcal{T}^{t}$, $i=1,\ldots, \ell$, and $\xi^{t,j} \in \mathcal{T}
^{t,j}$, $j \in a_{00}^- (\bar x) \backslash \{\bar j\}$, we obtain  $\beta_{\bar j}\xi^{t,\bar j}_{\bar j} = 0$. Due to (\ref{eq:Step5DF2ne0}), it must hold
$\beta_{\bar j}=0$.
For $\bar j \in a_{00}^{-,2}\left(\bar x\right)$ we get, by multiplying both sides with $e_{\kappa +\bar j}^T$ from the left, that $\beta_{\bar j}=0$, due to (\ref{eq:Step5DF1ne0}).
Overall, the second sum in the considered equation vanishes.
Hence, it simplifies to
\[
\sum\limits_{i=1}^{\ell} b_i \xi^{i} = 0.
\]
But then, due to the choice of $\xi^{i}$'s, we have $b_i=0$ for all $i=1,\ldots,\ell.$
Thus, the  claimed linear independence follows.

{\bf Step 5c.}
\begin{comment}
    In view of (\ref{eq:Step5DF2ne0}) and (\ref{eq:Step5DF1ne0}) we can assume
\[
DF_{1,j}\left(\bar x\right)\xi^{t,j}=1 \mbox{ for all }
j \in a_{00}^{-,1}\left(\bar x\right) \quad \mbox{and}\quad
DF_{2,j}\left(\bar x\right)\xi^{t,j}=1 \mbox{ for all }
j \in a_{00}^{-,2}\left(\bar x\right).
\]
\end{comment}
Suppose $\bar j \in a_{00}^{-}\left(\bar x\right)$ and $\xi^{t,\bar j}$ be chosen as before.
We claim that it holds:
\begin{equation}
\label{eq:step5c}
\lim\limits_{t \to 0} \left(\xi^{t,\bar j}\right)^T  D^2L^\mathcal{S}(x^t) \xi^{t,\bar j}= -\infty.
\end{equation}
To see this, let us recall from Step 2:
\[
\begin{array}{rcl}
\displaystyle
D^2 L^\mathcal{S}\left( x^t\right) &=& 
\displaystyle D^2 L\left( x^t\right) + \sum\limits_{j=1}^p\eta^t_j\left(e_j e_{\kappa+j}^T+e_{\kappa+j} e_{j}^T\right) \\&&+\displaystyle\sum\limits_{j=p+1}^\kappa \eta^t_j\left(e_j D F_{2,j}\left(x^t\right)+D^T F_{2,j}\left(x^t\right) e_j^T\right).
\end{array}
\]
By definition of $\xi^{t,\bar j}$, we conclude
\[
\begin{array}{rcl}
\displaystyle
\left(\xi^{t,\bar j}\right)^T D^2 L^\mathcal{S}\left( x^t\right) \xi^{t,\bar j} &=& 
\left(\xi^{t,\bar j}\right)^T \displaystyle D^2 L\left( x^t\right)\xi^{t,\bar j} 
+\left(\xi^{t,\bar j}\right)^T \left(\eta^t_{\bar j}\left(e_{\bar j} e_{\kappa+\bar j}^T+e_{\kappa+\bar j} e_{\bar j}^T\right) \right)\xi^{t,\bar j}.
\end{array}
\]
Note that for $t \to 0$ we have:
$$\left(\xi^{t,\bar j}\right)^T \displaystyle D^2 L\left( x^t\right)\xi^{t,\bar j}  \rightarrow \left(\xi^{\bar j}\right)^T \displaystyle D^2 L\left( \bar x\right)\xi^{\bar j}.$$
Hence, for our claim to be valid, it is sufficient to show that
\[
\lim\limits_{t \to 0}
\left(\xi^{t,\bar j}\right)^T \left(\eta^t_{\bar j}\left(e_{\bar j}\cdot e_{\kappa+\bar j}^T+e_{\kappa+\bar j}\cdot e_{\bar j}^T\right) \right)\xi^{t,\bar j}= -\infty.
\]
For that, let us first assume $\bar j \in a_{00}^{-,1}\left(\bar x\right)$. Due to $\xi^{t, \bar j} \in \mathcal{T}^{t,\bar j}$,
we can rewrite the term under consideration as follows
\[
\displaystyle -\eta^t_{\bar j}\left(\left(\xi^{t,\bar j}\right)^Te_{\bar j}\cdot \frac{x^t_{\bar \kappa +\bar j}}{x^t_{\bar j}}e_{\bar j}^T\xi^{t,\bar j}+\left(\xi^{t,\bar j}\right)^T\frac{x^t_{\bar \kappa +\bar j}}{x^t_{\bar j}}e_{\bar j}\cdot e_{\bar j}^T \xi^{t,\bar j}\right)
  =\displaystyle-2\eta^t_{\bar j}x^t_{\bar \kappa +\bar j}\frac{1}{x^t_{\bar j}} \left(\xi^{t,\bar j}\right)^T e_{\bar j}\cdot e_{\bar j}^T\xi^{t,\bar j}. 
  %\\  =&2\eta_{\bar j}F_{2,\bar j}\left(x^t\right) \frac{1}{F_{1,\bar j}\left(x^t\right)}.
\]
Since $\bar j \in a_{00}^-\left(\bar x\right)$, we have, due to Corollary \ref{cor:activeindex:scholtes} and Proposition \ref{prop:com:scholtes}, for $t \rightarrow 0$:
\[
 %- \eta^t_{\bar j} F_{2,\bar j}\left(x^t\right)=
 - \eta^t_{\bar j} x_{\kappa+\bar j}^t \rightarrow \bar \varrho_{1,\bar j} < 0,
\]
and, further,
\[
  \frac{1}{x_{\bar j}} \rightarrow +\infty.
\]
Finally, by construction
$$\left(\xi^{t,\bar j}\right)^T e_{\bar j}\cdot e_{\bar j}^T\xi^{t,\bar j} \rightarrow \left(\xi^{\bar j}_{\bar j}\right)^2 \not =0.$$ 
Overall, we conclude that assertion
(\ref{eq:step5c}) indeed holds for $\bar j \in a_{00}^{-,1}\left(\bar x\right)$.
If instead $\bar j \in a_{00}^{-,2}\left(\bar x\right)$, the proof is analogous.

{\bf Step 5d.}
We observe for the quadratic index of $x^t$ with $t$ sufficiently small
\[
\begin{array}{rcl}
q&=&QI^\mathcal{S}_{t,\mathcal{T}^\mathcal{S}_{x^t}}\overset{(\ref{eq:step3b})}{=} QI^\mathcal{S}_{t,\mathcal{T}^t} 
+QI^\mathcal{S}_{t,\mathcal{T}^\mathcal{S}_{x^t} \backslash \mathcal{T}^t}
\overset{(\text{Step 3a})}{=} QI^\mathcal{S}_{t,\mathcal{T}_{\bar x}} 
+QI^\mathcal{S}_{t,\mathcal{T}^\mathcal{S}_{x^t} \backslash \mathcal{T}^t}\\
&\overset{(\text{Steps 5a--c})}{\ge}&QI^\mathcal{S}_{t,\mathcal{T}_{\bar x}} 
+\left\vert a_{00}^-\left(\bar x\right)\right\vert
\overset{(\ref{eq:step2a})}{=}\overline{QI}_{\mathcal{T}_{\bar x}}
+\left\vert a_{00}^-\left(\bar x\right)\right\vert
=CI.
\end{array}
\]
%Note, that for $a_{10}\left(\bar x\right)\ne \emptyset$ we follow the exact same arguments to get (\ref{eq:main-index}).
%\end{proof}
\qed

\subsection*{Proof of Theorem \ref{thm:wellposedness}.}
For the sake of simplicity, we assume without loss of generality that 
$a_{10}\left(\bar x\right)=\emptyset$ and that the considered problems are given in the standard forms (\ref{eq:mpccstd}) and (\ref{eq:scholtesstd}), respectively.

(i) First, we show the existence of nondegenerate Karush-Kuhn-Tucker points of $\mathcal{S}$ in a neighborhood of $\bar x$.

{\bf Step 1.}
We consider the auxiliary system of equations $G(t,x,\sigma,\varrho) =0$ given by (\ref{eq:ift-stat})-(\ref{eq:ift-f4}), which mimics stationarity and feasibility. For stationarity we use:
\begin{equation}
\label{eq:ift-stat}
  \begin{array}{rl}
   -\nabla f(x) &  \displaystyle
 +\sum\limits_{j \in a_{01}^-\left(\bar x\right)}\frac{\sigma_{1,j}}{F_{2,j}\left(\bar x\right)}
\left(F_{2,j}\left(x\right) e_{j} + x_j D^T F_{2,j}\left(x\right)\right) +\sum\limits_{j \in a_{01}^+\left(\bar x\right)}\sigma_{1,j}
e_j\\ &
 %\displaystyle +\sum\limits_{j \in a_{10}^-\left(\bar x\right)} \frac{\sigma_{2,j}}{F_{1,j}\left(\bar x\right)}
%\left(F_{1,j}\left(x\right) \nabla F_{2,j}\left(x\right) + F_{2,j}\left(x\right) \nabla F_{1,j}\left(x\right)\right)  +\sum\limits_{j \in a_{10}^+\left(\bar x\right)}\sigma_{2,j}
%\nabla F_{2,j}\left( x\right)
%        \\ \\&
    +\displaystyle \sum\limits_{j \in a_{00}\left(\bar x\right)} \left( \varrho_{1,j}
e_j
   + \varrho_{2,j}e_{\kappa+j}\right) =0.
    \end{array}
\end{equation}
For feasibility we use:
\begin{equation}
\label{eq:ift-f1}
\frac{1}{F_{2,j}\left(\bar x\right)}\left(x_{j}\cdot F_{2,j}(x)-t\right)=0,\,j\in a_{01}^-\left(\bar x\right), \quad
    x_j=0,\,j\in a_{01}^+\left(\bar x\right),     
\end{equation}
%\begin{equation}
%\label{eq:ift-f2}
%\frac{1}{F_{1,j}\left(\bar x\right)}\left(F_{1,j}(x)\cdot F_{2,j}(x)-t\right)=0,\,j\in a_{10}^-\left(\bar x\right), \quad F_{2,j}\left(x\right)=0,\,j \in a_{10}^+\left(\bar x\right),   
%\end{equation}
\begin{equation}
\label{eq:ift-f3}
  x_j+\frac{\varrho_{2,j}\cdot \sqrt{t} }{\sqrt{\varrho_{1,j}\cdot \varrho_{2,j}}}=0,\,
   j \in a_{00}^-\left(\bar x\right), \quad
   x_j=0,\, 
    j \in a_{00}^+\left(\bar x\right),
\end{equation}
\begin{equation}
\label{eq:ift-f4}
  x_{\kappa+j}+\frac{\varrho_{1,j}\cdot \sqrt{t} }{\sqrt{\varrho_{1,j}\cdot \varrho_{2,j}}}=0,\,
      j \in a_{00}^-\left(\bar x\right), \quad 
     x_{\kappa+j}=0,\,
    j \in  a_{00}^+\left(\bar x\right).
\end{equation}
In view of feasibility and C-stationarity of $\bar x$ with multipliers $(\bar \sigma, \bar \rho)$, the vector
$(0,\bar x,\bar \sigma,\bar \varrho)$ solves the system of equations (\ref{eq:ift-stat})-(\ref{eq:ift-f4}).

{\bf Step 2.}
We consider the blockwise matrix
$$
\displaystyle \frac{\partial G (t,x, \sigma, \rho)}{\partial (x,\sigma,\varrho)}=\begin{bmatrix}
A&B\\
B^T&D
\end{bmatrix}.$$
Here, we have
\[
\begin{array}{rcl}
      A&=& \displaystyle - D^2f(x)%\\ \\ && 
      \displaystyle  
     
+ \begin{aligned}\sum\limits_{j \in a_{01}^-\left(\bar x\right)}
      \frac{\sigma_{1,j}}{F_{2,j}\left(\bar x\right)}\left(
D^T F_{2,j}\left(x\right) e_j^T+ 
 e_{j} DF_{2,j}\left(x\right) + 
 x_{j} D^2F_{2,j}\left(x\right)\right).
 \end{aligned}
% \\ \\&&
%  \displaystyle +\sum\limits_{j \in a_{01}^+\left(\bar x\right)}\sigma_{1,j}
%%%%%%D^2F_{1,j}\left( x\right)
%\\ \\&&
% \displaystyle
%       \begin{aligned}
%+ \sum\limits_{j \in a_{10}^-\left(\bar x\right)}
%      \frac{\sigma_{2,j}}{F_{1,j}\left(\bar x\right)}
%\biggl(&\nabla F_{2,j}\left(x\right) DF_{1,j}\left(x\right) + 
% F_{2,j}\left(x\right) D^2F_{1,j}\left(x\right)\\ &\quad+ 
% \nabla F_{1,j}\left(x\right) DF_{2,j}\left(x\right) + 
% F_{1,j}\left(x\right) D^2F_{2,j}\left(x\right) \biggr)
%  \end{aligned}
%\\ \\
% && \displaystyle + \sum\limits_{j \in a_{10}^+\left(\bar x\right)}\sigma_{2,j}
%D^2 F_{2,j}\left( x\right)
%%%%%+ \displaystyle \sum\limits_{j \in a_{00}\left(\bar x\right)} \left( \varrho_{1,j}
%D^2F_{1,j}\left(x\right)
%   + \varrho_{2,j}D^2F_{2,j}\left(x\right)\right).
\end{array}
\]
The columns of $B$ are given by the vectors:
\[
\frac{1}{F_{2,j}\left(\bar x\right)}
\left(F_{2,j}\left(x\right) e_{j} + x_{j} D^T F_{2,j}\left(x\right)\right), j \in a_{01}^-\left(\bar x\right),\quad
e_{j}, j \in a_{01}^+\left(\bar x\right),
\quad
%\[
%\frac{1}{F_{1,j}\left(\bar x\right)}
%\left(F_{1,j}\left(x\right) \nabla F_{2,j}\left(x\right) + F_{2,j}\left(x\right) \nabla F_{1,j}\left(x\right)\right), j \in a_{10}^-\left(\bar x\right), \quad
%\nabla F_{2,j}\left( x\right), j \in a_{10}^+\left(\bar x\right),
%\]
e_{j}, e_{\kappa+j}, j \in a_{00}\left(\bar x\right).
\]
The first $\left\vert a_{01}\left(\bar x\right)\right\vert%+\left\vert a_{10}\left(\bar x\right)\right\vert
$ rows of $D$ are vanishing. The subsequent rows of $D$ are given by the vectors:
\[
\begin{pmatrix}
0\\
    -\frac{\sqrt{-\varrho_{2,j}t}}{2\sqrt{-\varrho_{1,j}^3}}e_j\\
    \frac{\sqrt{t}}{2\sqrt{\varrho_{1,j}\varrho_{2,j}}}e_j
\end{pmatrix},\,j \in a_{00}^-\left(\bar x\right),
\]
which are followed by $\left|a_{00}^+\left(\bar x\right)\right|$ vanishing rows.
Finally, the next rows of $D$  are given by the vectors:
\[
\begin{pmatrix}
0\\
    \frac{\sqrt{t}}{2\sqrt{\varrho_{1,j}\varrho_{2,j}}}e_j\\
    -\frac{\sqrt{-\varrho_{1,j}t}}{2\sqrt{-\varrho_{2,j}^3}}e_j
\end{pmatrix},\,j \in a_{00}^-\left(\bar x\right),
\]
while the remaining $\left|a_{00}^+\left(\bar x\right)\right|$ rows are vanishing.
Additionally we have $D=0$ at $(0,\bar x,\bar \sigma,\bar \varrho)$. Hence, we can apply Theorem 2.3.2 from \cite{jongen:2004}, which says that
$
%D G (0,\bar x,\bar \sigma, \bar \rho)=
\begin{bmatrix}
A&B\\
B^T&0
\end{bmatrix}$ is nonsingular if and only if
$\xi^TA\xi\ne0$ for all $\xi \in B^\perp\backslash \{0\}$. Here, $B^{\perp}$ refers to the orthogonal complement of the subspace spanned by the columns of $B$. 
In view of $B^\perp=\mathcal{T}_{\bar x}$ at $\left(0, \bar x, \bar \sigma, \bar \varrho\right)$, 
we check for $\xi \in B^{\perp}$ with $ \xi \ne 0$:
\[
\xi^TA\xi = \xi^T D^2 L\left(\bar x\right) \xi \overset{\text{NDC3}}{\ne} 0.
\]
Hence, by means of the implicit function theorem we obtain for any $t>0$ sufficiently small a solution $\left(t,x^t,\sigma^t,\varrho^t\right)$ of the system of equations (\ref{eq:ift-stat})-(\ref{eq:ift-f4}).

{\bf Step 3.}
By choosing $t$ even smaller if necessary, we can ensure due to continuity reasons as well as  NDC2 and NDC4 that the following holds:
\begin{itemize}
    \item[(i)] $\mbox{sgn}\left(\sigma^t_{1,j}\right)=\mbox{sgn}\left(\bar \sigma_{1,j}\right)$, $F_{2,j}\left(x^t\right)>0$, \,$j\in a_{01}\left(\bar x\right)$,
   \item[(ii)] %$\mbox{sgn}\left(\sigma^t_{2,j}\right)=\mbox{sgn}\left(\bar \sigma_{2,j}\right)$, $F_{1,j}\left(x^t\right)>0$,\,$j \in a_{10}\left(\bar x\right)$,
   %\item[(iii)] 
   $\mbox{sgn}\left(\varrho^t_{1,j}\right)=\mbox{sgn}\left(\bar \varrho_{1,j}\right)$, $ \mbox{sgn}\left(\varrho^t_{2,j}\right)=\mbox{sgn}\left(\bar \varrho_{2,j}\right)$,\,$j\in a_{00}\left(\bar x\right)$.
\end{itemize} 
From here it is straightforward to see that $x^t \in M^\mathcal{S}$ is feasible for $\mathcal{S}$ and we have:
\begin{itemize}
    \item[(a)] $\mathcal{H}\left(x^t\right)=a_{01}^-\left(\bar x\right)%\cup a_{10}^-\left(\bar x\right)
    \cup a_{00}^-\left(\bar x\right)$,
    \item[(b)] $\mathcal{N}_1\left(x^t\right)=a_{01}^+\left(\bar x\right)\cup a_{00}^+\left(\bar x\right)$,
    \item[(c)] $\mathcal{N}_2\left(x^t\right)=%a_{10}^+\left(\bar x\right)\cup
    a_{00}^+\left(\bar x\right)$.    
\end{itemize}

Further, it holds:
\begin{equation}
    \label{eq:help34}
   \begin{array}{rcl}
   \nabla f\left(x^t\right) &=&\displaystyle 
 \sum\limits_{j \in a_{01}^-\left(\bar x\right)}\frac{\sigma^t_{1,j}}{F_{2,j}\left(\bar x\right)}
\left(F_{2,j}\left(x^t\right) e_{j} + x_{j} D^T F_{2,j}\left(x^t\right)\right) 
%\\ \\&&
 %\displaystyle +\sum\limits_{j \in a_{10}^-\left(\bar x\right)} \frac{\sigma^t_{2,j}}{F_{1,j}\left(\bar x\right)}
%\left(F_{1,j}\left(x^t\right) \nabla F_{2,j}\left(x^t\right) + F_{2,j}\left(x^t\right) \nabla F_{1,j}\left(x^t\right)\right)  
+\sum\limits_{j \in a_{01}^+\left(\bar x\right)}\sigma^t_{1,j}
e_{j}\\ &&
\displaystyle
%+\sum\limits_{j \in a_{10}^+\left(\bar x\right)}\sigma^t_{2,j}
%\nabla F_{2,j}\left( x^t\right)
        %\\ \\&&
    +\displaystyle \sum\limits_{j \in a_{00}\left(\bar x\right)} \left( \varrho_{1,j}^t
e_{1,j}
   + \varrho_{2,j}^t e_{\kappa+j}\right).
    \end{array}
\end{equation}
We rename the multipliers as follows:
\[
\eta^{t}_j=\left\{
\begin{array}{ll}
\displaystyle -\frac{\sigma^t_{1,j}}{F_{2,j}\left(\bar x\right)}& \mbox{for }j\in a_{01}^-\left(\bar x\right),\\ 
%\displaystyle -\frac{\sigma^t_{2,j}}{F_{1,j}\left(\bar x\right)},& \mbox{for }j\in a_{10}^-\left(\bar x\right),\\ 
\displaystyle -\frac{\varrho^t_{1,j}x^t_{j}}{t}& \mbox{for }j\in a_{00}^-\left(\bar x\right),\\
0&\mbox{else,}
\end{array}\right.
\]
\[
\nu^t_{1,j}=\left\{
\begin{array}{ll}
\displaystyle \sigma_{1,j}^t& \mbox{for }j\in a_{01}^+\left(\bar x\right),\\ 
\displaystyle \varrho_{1,j}^t& \mbox{for }j\in a_{00}^+\left(\bar x\right),\\ 
0&\mbox{else,}
\end{array}\right.\\
\quad
\nu^t_{2,j}=\left\{
\begin{array}{ll}
%\displaystyle \sigma_{2,j}^t,& \mbox{for }j\in a_{10}^+\left(\bar x\right),\\ 
\displaystyle \varrho_{2,j}^t& \mbox{for }j\in a_{00}^+\left(\bar x\right),\\ 
0&\mbox{else.}
\end{array}\right.\\
\]
By doing so, it is straightforward to confirm that for all $ j \in  a_{00}^-\left(\bar x\right)$ it then holds:
\[
%\frac{\varrho^t_{1,j}F_{1,j}\left(x^t\right)}{t}=\frac{\varrho^t_{2,j}F_{2,j}\left(x^t\right)}{t}
\varrho^t_{1,j}=-\eta^t_j x^t_{\kappa+j}\quad\text{and}\quad \varrho^t_{2,j}=-\eta^t_j x^t_{j}. 
\]
Substituting into (\ref{eq:help34}), we conclude that $x^t$ fulfills (\ref{eq:kkt-1}) and (\ref{eq:kkt-2}) with multipliers $(\eta^t,\nu^t)$. 
Consequently, $x^t$ is a Karush-Kuhn-Tucker point of $\mathcal{S}$. 

{\bf Step 4.}
In view of Lemma \ref{lem:LICQ-MPCCLICQ}, ND1 is satisfied at any $x^t$ for $t$ sufficiently small. Obviously, ND2 holds as well.
Therefore, it remains to show ND3, i.e.~the restriction of 
$D^2 L^\mathcal{S}\left( x^t\right)$ on $\mathcal{T}^{\mathcal{S}}_{x^t}$ is nonsingular. Note that $\mathcal{T}^{\mathcal{S}}_{x^t}$ is of dimension $n-\alpha_t^{\mathcal{S}}$, cf. Step 1 from the proof of Theorem \ref{thm:kktsequence}. In view of (a)-(c), we get from Step 3 here that $\alpha_t^{\mathcal{S}}$ is constant for $t$ sufficiently small. Thus, we refer to it as $\alpha^{\mathcal{S}}$.
Next, we construct a basis of $\mathcal{T}^{\mathcal{S}}_{x^t}$ as follows. 
First, we choose eigenvectors $\bar \xi_1,\ldots,\bar \xi_{n-\alpha}$ of $D^2 L\left(\bar x\right)$ forming a basis of $\mathcal{T}_{\bar x}$, cf. Step 1 from the proof of Theorem \ref{thm:kktsequence}.
With similar arguments as in Step 3 of the proof of Theorem \ref{thm:kktsequence} we find $\xi_k^t \in \mathcal{T}^{\mathcal{S}}_{x^t}$, $k=1,\ldots,n-\alpha$, still linearly independent. 
For the remaining $n-\alpha^{\mathcal{S}}-\left(n-\alpha\right)=\left\vert a^-_{00}\left(\bar x\right)\right\vert$ vectors we chose $\xi^{t,\bar j}, \bar j \in  a_{00}^-(\bar x)$ as constructed in Step 5a of the proof of Theorem \ref{thm:kktsequence}. Due to the same arguments as in Step 5b there, the chosen vectors do indeed form a basis of $\mathcal{T}^{\mathcal{S}}_{x^t}$ together.

Further, due to (\ref{eq:step2limits}) and NDC3, we have for $k=1,\ldots,n-\alpha$:
\[\text{sgn}\left(\left(\bar \xi_k\right)^TD^2L^{\mathcal{S}}\left(x^t\right)\bar \xi_k\right)
=\text{sgn} \left(\left(\bar \xi_k\right)^TD^2L\left(\bar x\right)\bar \xi_k\right)
\ne 0.
\]
Hence, by construction it also holds for $t$ sufficiently small:
\[
\text{sgn}\left(\left(\xi_k^t\right)^TD^2L^{\mathcal{S}}\left(x^t\right)\xi_k^t\right)
=\text{sgn}\left(\left(\bar \xi_k\right)^TD^2L\left(\bar x\right)\bar \xi_k\right)\ne 0.
\]
We argue in the same way as in Step 5c of the proof of Theorem \ref{thm:kktsequence}
to conclude that for $t$ sufficiently small and all $\bar j \in a_{00}^-(\bar x)$ it holds:
\[
\left(\xi^{t,\bar j}\right)^TD^2L^{\mathcal{S}}\left(x^t\right)\xi^{t,\bar j}<0.
\]
Combining the latter results yields that ND3 is fulfilled at $x^t$. Moreover, it also follows directly for the corresponding quadratic index that
\[
QI^\mathcal{S}_{t,\mathcal{T}^\mathcal{S}_{x^t}} =\overline{QI}_{t,\mathcal{T}_{\bar x}}+\left\vert a_{00}^-\left(\bar x\right)\right\vert = c.
\]

(ii) Next, we elaborate on the uniqueness of Karush-Kuhn-Tucker points $x^t$ constructed above. 

{\bf Step 5.}
Let us consider any Karush-Kuhn-Tucker points $\widetilde x^t \in M^\mathcal{S}$ of $\mathcal{S}$ with multipliers $(\widetilde \eta^t,\widetilde \nu^t)$ in a sufficiently small neighborhood of $\bar x$. In view of Proposition \ref{prop:com:scholtes}, NDC2, and NDC4, the relations in Corollary \ref{cor:activeindex:scholtes} --
if applied to $\widetilde x^t$ -- provide:
\[a_{01}^-\left(\bar x\right) \cup a_{00}^-\left(\bar x\right)= \mathcal{H}\left(\widetilde x^t\right), \quad a_{01}^+\left(\bar x\right)= \mathcal{N}_1\left(\widetilde x^t\right) \backslash \mathcal{N}_2\left(\widetilde x^t\right), \quad
a_{00}^+\left(\bar x\right)= \mathcal{N}_1\left(\widetilde x^t\right) \cap \mathcal{N}_2\left(\widetilde x^t\right).
\]

Further, $\widetilde x^t$ together with $(\widetilde \eta^t,\widetilde \nu^t)$ has to fulfill (\ref{eq:kkt-1}) and (\ref{eq:kkt-2}).
We define the multipliers
\[
\begin{array}{c}
\widetilde \sigma_{1,j}^t=-\widetilde \eta_j  F_{2,j}\left(\bar x\right), j \in a_{01}^-(\bar x),  \quad
\widetilde \sigma_{1,j}^t=\widetilde\nu_{1,j}j \in a_{01}^+(\bar x),\\
%\widetilde \sigma_{2,j}^t&=&-\widetilde \eta_j  F_{1,j}\left(\bar x\right),&j \in a_{10}^-(\bar x),\\
%\widetilde \sigma_{2,j}^t&=&\widetilde\nu_{2,j},&j \in a_{10}^+(\bar x),\\
\widetilde \varrho_{1,j}^t=-\widetilde \eta_j  \bar x_{\kappa+j}, j \in a_{00}^-(\bar x), \quad 
\widetilde \varrho_{1,j}^t=\widetilde\nu_{1,j},j \in a_{00}^+(\bar x),\\
\widetilde \varrho_{2,j}^t=-\widetilde \eta_j  \bar x_{j},j \in a_{00}^-(\bar x), \quad
\widetilde \varrho_{2,j}^t=\widetilde\nu_{2,j}, j \in a_{00}^+(\bar x).\\
\end{array}
\]
A straightforward calculation shows that $\left(t, \widetilde x^t, \widetilde \sigma^t, \widetilde \varrho^t\right)$ fulfills equations (\ref{eq:ift-stat})-(\ref{eq:ift-f4}) for $t$ sufficiently small. However, the implicit function theorem was used in Step 2 to find the solution of this system of equations in the neighborhood of $\bar x$. Hence, $x^t$ must be unique.    \qed

\color{black}

\section*{Acknowledgements}
We would like to express our sincere gratitude to the anonymous referees who helped to considerably improve the quality of the article and to properly contextualize the new and the existing results. 
}

%\THEEndNotes
%\begingroup \parindent 0pt \parskip 4ex
%\def\enotesize{\normalsize} 
%\theendnotes
%\endgroup

% Appendix here
% Options are (1) APPENDIX (with or without general title) or
%             (2) APPENDICES (if it has more than one unrelated sections)
% Outcomment the appropriate case if necessary
%
% \begin{APPENDIX}{<Title of the Appendix>}
% \end{APPENDIX}
%
%   or
%
% \begin{APPENDICES}
% \section{<Title of Section A>}
% \section{<Title of Section B>}
% etc
% \end{APPENDICES}

% Acknowledgments here
%\ACKNOWLEDGMENT{We would like to express our sincere gratitude to [acknowledge individuals, organizations, or institutions] for their invaluable contributions to this research. We are also grateful to [mention any additional acknowledgements, such as technical assistance, data providers, or colleagues] for their support and assistance throughout the course of this work.}

% References here (outcomment the appropriate case)

% CASE 1: BiBTeX used to constantly update the references
%   (while the paper is being written).
%\bibliographystyle{informs2014} % outcomment this and next line in Case 1
%\bibliography{<your bib file(s)>} % if more than one, comma separated

%%%%%%%%%%%%%%%%%%%%%%%%%%%%%%%%%%%%%%
%References
%%%%%%%%%%%%%%%%%%%%%%%%%%%%%%%%%%%%%%
\bibliographystyle{apalike}
\bibliography{lit}

% CASE 2: BiBTeX used to generate mypaper.bbl (to be further fine tuned)
%\input{mypaper.bbl} % outcomment this line in Case 2

%If you don't use BiBTex, you can manually itemize references as shown below.

%\bibliographystyle{nonumber}

%%%%%%%%%%%%%%%%%
\end{document}